\documentclass[reqno]{amsart}

\usepackage{amsthm}

\usepackage{
amssymb,amsmath,enumerate,stackrel,tabu,mathrsfs,enumitem,verbatim,stmaryrd,
leftidx}
\usepackage{xcolor} 
\usepackage{xr-hyper}
\usepackage{nicefrac,stackrel}
\usepackage[all,cmtip]{xy}
\usepackage[utf8]{inputenc} 
\chardef\cprime"7E 
\usepackage{cases}

\usepackage[T1]{fontenc} 

\usepackage[pagebackref,hyperindex,linktocpage=true]{hyperref}
\hypersetup{
    colorlinks,
    linkcolor={red!50!black},
    citecolor={red!50!black}, 
    urlcolor={red!50!black}, 
    filecolor={red!50!black} 
}

\usepackage{tikz}

\usepackage[top=1in, left=1in, right=1in, bottom=1in, includeheadfoot]{geometry}

\usepackage{xcolor}

\def\label#1{\label{#1}}
\def\mylabel#1{\thlabel{#1}}
\def\myref#1{\thref{#1}}

\usepackage{theoremref}
\definecolor{labelkey}{rgb}{1,0,0}

\setlist[enumerate,1]{label=\roman*\textup{)},
style=nextline, leftmargin=7mm,
topsep=0mm,
itemsep=0.4mm, 
align=parleft
}
\makeatletter
\@addtoreset{equation}{section}
\makeatother

\numberwithin{equation}{section}

\theoremstyle{definition}
\newtheorem{Defi}{Definition}[section] \newcommand{\defi}{\begin{Defi}} \newcommand{\xdefi}{\end{Defi}} 
\newtheorem{Cons}[Defi]{Construction} \newcommand{\cons}{\begin{Cons}} \newcommand{\xcons}{\end{Cons}} 
\newtheorem{DefiLemm}[Defi]{Definition and Lemma} \newcommand{\defilemm}{\begin{DefiLemm}} \newcommand{\xdefilemm}{\end{DefiLemm}} 
\newtheorem{Bsp}[Defi]{Example} \newcommand{\exam}{\begin{Bsp}} \newcommand{\xexam}{\end{Bsp}} 
\newtheorem{Syno}[Defi]{Synopsis} \newcommand{\syno}{\begin{Syno}} \newcommand{\xsyno}{\end{Syno}} 
\newtheorem{Bem}[Defi]{Remark} \newcommand{\rema}{\begin{Bem}} \newcommand{\xrema}{\end{Bem}} 

\theoremstyle{plain}
\newtheorem{Theo}[Defi]{Theorem} \newcommand{\theo}{\begin{Theo}} \newcommand{\xtheo}{\end{Theo}} 
\newtheorem{Satz}[Defi]{Proposition} \newcommand{\prop}{\begin{Satz}} \newcommand{\xprop}{\end{Satz}} 
\newtheorem{Lemm}[Defi]{Lemma} \newcommand{\lemm}{\begin{Lemm}} \newcommand{\xlemm}{\end{Lemm}} 
\newtheorem{Coro}[Defi]{Corollary} \newcommand{\coro}{\begin{Coro}} \newcommand{\xcoro}{\end{Coro}} 

\newcommand{\refsect}[1]{§\ref{sect--#1}}
\newcommand{\refit}[1]{\ref{item--#1}}
\newcommand{\refeq}[1]{(\ref{eqn--#1})}

\newcommand{\eqn}{\begin{equation}} \newcommand{\xeqn}{\end{equation}}
\newcommand{\eqnarr}{\begin{eqnarray*}} \newcommand{\xeqnarr}{\end{eqnarray*}}
\newcommand{\eqnarra}{\begin{eqnarray}} \newcommand{\xeqnarra}{\end{eqnarray}}

\newcommand{\pf}{\begin{proof}} \newcommand{\xpf}{\end{proof}}

\newcommand{\nc}{\newcommand}
\nc{\StP}[1]{\cite[\href{http://stacks.math.columbia.edu/tag/#1}{Tag #1}]{StacksProject}}\nc{\StPd}[3]{\cite[Tags~\href{http://stacks.math.columbia.edu/tag/#1}{#1}, 
\href{http://stacks.math.columbia.edu/tag/#2}{#2}, 
\href{http://stacks.math.columbia.edu/tag/#3}{#3}]{StacksProject}} 

\nc{\on}{\operatorname}
\nc{\aff}{{\on{aff}}}
\nc{\modi}{{\on{mod}}} 
\nc{\even}{{\on{even}}}
\nc{\odd}{{\on{odd}}}
\nc{\naive}{{\on{naive}}}
\nc{\hofib}{\on{hofib}}
\nc{\Bun}{\on{Bun}}
\nc{\ad}{{\on{ad}}}
\nc{\lft}{{\on{lft}}}

\nc{\str}{\on{-}}
\nc{\perf}{{\on{perf}}}
\nc{\Rel}{{\on{Pos}}}
\nc{\lan}{\langle}
\nc{\ran}{\rangle}

\nc{\bbA}{{\A}} 
\nc{\bbB}{{\mathbb B}}
\nc{\bbC}{{\mathbb C}}
\nc{\bbD}{{\mathbb D}}
\nc{\bbE}{{\mathbb E}}
\nc{\bbG}{{\mathbf G}}
\nc{\bbH}{{\mathbb H}}
\nc{\bbI}{{\mathbb I}}
\nc{\bbJ}{{\mathbb J}}
\nc{\bbK}{{\mathbb K}}
\nc{\bbL}{{\mathbb L}}
\nc{\bbM}{{\mathbb M}}
\nc{\bbN}{{\N}} 
\nc{\bbO}{{\mathbb O}}
\nc{\bbP}{{\P}} 
\nc{\bbQ}{{\Q}} 
\nc{\bbR}{{\mathbb R}}
\nc{\bbS}{{\mathbb S}}
\nc{\bbT}{{\mathbb T}}
\nc{\bbU}{{\mathbb U}}
\nc{\bbV}{{\mathbb V}}
\nc{\bbW}{{\mathbb W}}
\nc{\bbX}{{\mathbb X}}
\nc{\bbY}{{\mathbb Y}}
\nc{\bbZ}{{\bf Z}}

\nc{\calA}{{\mathcal A}}
\nc{\calB}{{\mathcal B}}
\nc{\calC}{{\mathcal C}}
\nc{\calD}{{\mathcal D}}
\nc{\calE}{{\mathcal E}}
\nc{\calF}{{\mathcal F}}
\nc{\calG}{{\mathcal G}}
\nc{\calH}{{\mathcal H}}
\nc{\calI}{{\mathcal I}}
\nc{\calJ}{{\mathcal J}}
\nc{\calK}{{\mathcal K}}
\nc{\calL}{{\mathcal L}}
\nc{\calM}{{\mathcal M}}
\nc{\calN}{{\mathcal N}}
\nc{\calO}{{\mathcal O}}
\nc{\calP}{{\mathcal P}}
\nc{\calQ}{{\mathcal Q}}
\nc{\calR}{{\mathcal R}}
\nc{\calS}{{\mathcal S}}
\nc{\calT}{{\mathcal T}}
\nc{\calU}{{\mathcal U}}
\nc{\calV}{{\mathcal V}}
\nc{\calW}{{\mathcal W}}
\nc{\calX}{{\mathcal X}}
\nc{\calY}{{\mathcal Y}}
\nc{\calZ}{{\mathcal Z}}

\nc{\Sht}{{\on{Sht}}}
\nc{\Frob}{{\on{Frob}}}
\nc{\Hecke}{{\on{Hecke}}}
\nc{\inv}{{\on{inv}}}
\nc{\Conv}{{\on{Conv}}}
\nc{\triv}{{\on{triv}}}
\nc{\Isom}{{\on{Isom}}}

\nc{\scrB}{{\mathscr{B}}}
\nc{\scrA}{{\mathscr{A}}}
\nc{\bbf}{{\mathbf{f}}}
\nc{\bba}{{\mathbf{a}}}

\nc{\al}{\alpha}
\nc{\be}{\beta}
\nc{\ga}{\gamma}
\nc{\la}{\lambda}
\nc{\qcqs}{{\on{qcqs}}}
\nc{\Bmu}{{\boldsymbol \mu}}

\nc{\pot}[1]{ [\hspace{-0,5mm}[ {#1} ]\hspace{-0,5mm}] }
\nc{\rpot}[1]{ (\hspace{-0,7mm}( {#1} )\hspace{-0,7mm}) }

\nc{\defined}{\hspace{0.1cm}\stackrel{\text{\tiny \rm def}}{=}\hspace{0.1cm}}
\nc{\co}{\colon}

\newcommand{\category}[1]{\mathrm{#1}}
\newcommand{\DGCat}{\category{DGCat}} 
\newcommand{\cont}{\category{cont}} 
\newcommand{\Ho}{\category{Ho}} 
\newcommand{\res}{\mathrm{res}} 
\newcommand{\Fun}{\category{Fun}} 
\newcommand{\AffSch}{\category{AffSch}} 
\newcommand{\PreStk}{\category{PreStk}} 
\newcommand{\Gpd}{\category{Gpd}} 
\newcommand{\Ani}{\category{Ani}} 
\newcommand{\Zar}{\mathrm{Zar}} 
\newcommand{\Du}{\mathrm{D}} 
\newcommand{\Sch}{\category{Sch}} 
\newcommand{\IndSch}{\category{IndSch}} 
\newcommand{\W}{\mathrm {W}} 

\def\wt{\mathrm {w}} 
\def\perv{\mathrm {perv}} 
\def\gr{\mathrm {gr}} 
\def\pH{^{\mathrm p} \H} 
\def\cl{\mathrm {cl}} 
\def\mot{\mathrm {m}} 
\def\motH{{}^\mot \H} 
\def\Gm{\mathbf {G}_\mathrm m} 
\newcommand{\GaX}[  1]{\mathbf {G}_{\mathrm {a}, #1}} 
\def\SL{\mathrm {SL}} 
\def\PGL{\mathrm {PGL}} 
\def\GL{\mathrm {GL}} 
\newcommand{\GmX}[  1]{\mathbf {G}_{\mathrm {m}, #1}} 
\def\IC{\mathrm{IC}} 
\def\red{\mathrm{red}} 
\def\pfp{\mathrm{pfp}} 
\def\ft{\mathrm{ft}} 

\font\tencyr=wncyr10
\font\sevencyr=wncyr7
\font\fivecyr=wncyr5
\newcommand{\colim}{\operatornamewithlimits{colim}} 
\def\id{{\rm id}} 
\def\pr{{\rm pr}} 
\def\opp{{\rm op}} 
\def\To#1#2{\mathop{\count0=#1 \loop\ifnum\count0>0 \smash-\mkern-7mu \advance\count0 -1 \repeat \mathord\rightarrow}\limits^{#2}} 
\def\Char{\mathop{\rm char}\nolimits} 
\def\CH{\mathop{\rm CH}\nolimits} 
\def\laxlim{\mathop{\rm laxlim}} 
\def\Hom{\mathop{\rm Hom}\nolimits} 
\def\bsl{\backslash} 
\def\Ind{\category{Ind}} 
\def\Gr{\mathop{\rm Gr}\nolimits} 
\def\Fl{\mathop{\rm Fl}\nolimits} 
\def\Sht{\mathop{\rm Sht}\nolimits} 
\def\Ext{\mathop{\rm Ext}\nolimits} 
\def\IHom{\underline{\Hom}} 
\def\Aut{\mathop{\rm Aut}\nolimits} 
\def\Rep{\category{Rep}} 
\def\bound{{\rm b}} 
\def\et{\mathrm{\acute et}} 
\def\incl{\mathrm{incl}} 

\definecolor{hellgrau}{RGB}{200,200,200} 
\definecolor{dunkelgrau}{RGB}{160,160,160} 
\definecolor{hellblau}{RGB}{194, 215, 249} %
\definecolor{dunkelblau}{RGB}{68, 128, 226} %

\def\Z{{\bf Z}} 
\def\bbF{{\bf F}}
\def\Fp{{\bf F}_p} %
\def\Fq{{\bf F}_q} %
\def\N{{\bf N}} 
\def\Q{{\bf Q}} 
\def\Qp{\Q_p} 
\def\Zp{\Z_p} 
\def\Ql{{\Q_\ell}} 
\def\Qq{{\overline \Q}} 
\def\A{{\bf A}} 
\renewcommand{\P}[1][1]{\mathbf P^{#1}} 
\def\Gm{\mathbf {G}_\mathrm m} 

\def\H{{\rm H}} 
\def\im{{\rm im}} 
\def\SH{\category{SH}} %
\def\DM{\category{DM}} 
\def\DTM{\category{DTM}} 
\def\Vect{\category{Vect}} 
\def\Perv{\category{Perv}} 
\def\MTM{\category{MTM}} 
\def\ii{$\infty$}
\def\bound{{\rm b}} 

\def\Spec{\mathop{\rm Spec}} 
\newcommand{\comp}{\mathrm{c}} 
\newcommand{\cstr}{\mathrm{cons}} 
\newcommand{\D}{\category{D}} 

\def\R{{\rm R}} 
\def\sbuildrel#1\over#2{\mathrel{\smash{\mathop{\kern0pt #2}\limits^{#1}}}}

\let\x\times

\renewcommand{\t}{\otimes}

\newcommand{\xtw}{\widetilde \x}
\renewcommand{\r}{\rightarrow}
\newcommand{\lr}{\longrightarrow}

\def\matrix#1{\null\,\vcenter{\normalbaselines
    \ialign{\hfil$##$\hfil&&\quad\hfil$##$\hfil\crcr
      \mathstrut\crcr\noalign{\kern-\baselineskip}
      #1\crcr\mathstrut\crcr\noalign{\kern-\baselineskip}}}\,}

\newdimen\harrowsize
\def\mapright#1{\smash{\mathop{\hbox to\harrowsize{\rightarrowfill}}\limits^{#1}}}

{\catcode`@=11
\gdef\cal{\fam\tw@}
\global\let\over\@@over
\global\let\atop\@@atop
\global\let\above\@@above
\global\let\overwithdelims\@@overwithdelims
\global\let\atopwithdelims\@@atopwithdelims
\global\let\abovewithdelims\@@abovewithdelims
\gdef\eqalign#1{\null\,\vcenter{\openup\jot\m@th
  \ialign{\strut\hfil$\displaystyle{##}$&$\displaystyle{{}##}$\hfil
      \crcr#1\crcr}}\,}
\newskip\xcentering \global\xcentering=0pt plus 1000pt minus 1000pt
\gdef\eqalignno#1{\displ@y \tabskip\xcentering
  \halign to\displaywidth{\hfil$\@lign\displaystyle{##}$\tabskip\z@skip
    &$\@lign\displaystyle{{}##}$\hfil\tabskip\xcentering
    &\llap{$\@lign##$}\tabskip\z@skip\crcr
    #1\crcr}}
\global\def\cases#1{\left\{\,\vcenter{\normalbaselines\m@th
    \ialign{$##\hfil$&\quad##\hfil\crcr#1\crcr}}\right.}
\gdef\eqlabel#1{\refstepcounter{equation}\label{eqn--#1}\eqno\hbox{\@eqnnum}}
}

\def \nts#1{ 
}

\def \journal#1
{ 
\noindent\colorbox{dunkelblau}{\parbox{\dimexpr\textwidth-2\fboxsep\relax}{#1}}
}

\usepackage{ifthen}

\def\iref#1{\ifthenelse{\equal{#1}{Bar_Construction}}{3.2.17}{}\ifthenelse{\equal{#1}{Chevalley_Triple}}{4.1.1}{}\ifthenelse{\equal{#1}{coro--DTM.G.X.generators}}{3.2.24}{}\ifthenelse{\equal{#1}{coro--equivalence.torsors}}{2.2.24}{}\ifthenelse{\equal{#1}{coro--equivariant.IC}}{5.3.6}{}\ifthenelse{\equal{#1}{coro--etale.torsor}}{A.4.8}{}\ifthenelse{\equal{#1}{coro--exactness}}{3.2.10}{}\ifthenelse{\equal{#1}{coro--Fl.stratification}}{4.3.12}{}\ifthenelse{\equal{#1}{coro--Ind.Artin.co.limit}}{2.3.4}{}\ifthenelse{\equal{#1}{coro--MTM.G.trivial}}{3.2.21}{}\ifthenelse{\equal{#1}{coro--prelim.intersection.motives}}{6.3.5}{}\ifthenelse{\equal{#1}{coro--t-structure}}{3.2.6}{}\ifthenelse{\equal{#1}{decomp_element}}{4.2.14}{}\ifthenelse{\equal{#1}{defi--adm.act}}{A.3.1}{}\ifthenelse{\equal{#1}{defi--cellular}}{3.1.5}{}\ifthenelse{\equal{#1}{defi--DM.G}}{2.2.6}{}\ifthenelse{\equal{#1}{defi--DM.prestacks}}{2.2.1}{}\ifthenelse{\equal{#1}{defi--DTM.G}}{3.1.21}{}\ifthenelse{\equal{#1}{defi--flag.variety}}{4.3.1}{}\ifthenelse{\equal{#1}{defi--intersection.motive}}{6.3.4}{}\ifthenelse{\equal{#1}{defi--motive.ind.scheme}}{2.3.12}{}\ifthenelse{\equal{#1}{defi--MTM.G}}{3.2.14}{}\ifthenelse{\equal{#1}{defi--rel.pos}}{6.1.2}{}\ifthenelse{\equal{#1}{defi--Schubert.scheme}}{4.4.1}{}\ifthenelse{\equal{#1}{defi--Schubert}}{4.3.4}{}\ifthenelse{\equal{#1}{defi--stratified.dfn}}{3.1.1}{}\ifthenelse{\equal{#1}{defi--stratified.G.action}}{3.1.26}{}\ifthenelse{\equal{#1}{defi--Tate.geometry}}{3.1.8}{}\ifthenelse{\equal{#1}{defi--unipotent}}{A.4.5}{}\ifthenelse{\equal{#1}{defi--Whitney-Tate.map}}{3.1.15}{}\ifthenelse{\equal{#1}{defilemm--Whitney.Tate.condition}}{3.1.11}{}\ifthenelse{\equal{#1}{deformation}}{A.4.10}{}\ifthenelse{\equal{#1}{dom_weights}}{4.1.2}{}\ifthenelse{\equal{#1}{double_classes}}{4.2.15}{}\ifthenelse{\equal{#1}{double_quotient}}{5.3.1}{}\ifthenelse{\equal{#1}{DTM.flag}}{5.2.1}{}\ifthenelse{\equal{#1}{eqn--adjunction.shriek}}{2.1.3}{}\ifthenelse{\equal{#1}{eqn--adjunction.star}}{2.1.2}{}\ifthenelse{\equal{#1}{eqn--Bar.descent}}{2.2.8}{}\ifthenelse{\equal{#1}{eqn--base.change.1}}{2.1.6}{}\ifthenelse{\equal{#1}{eqn--base.change.2}}{2.1.7}{}\ifthenelse{\equal{#1}{eqn--Beck}}{2.1.12}{}\ifthenelse{\equal{#1}{eqn--BS.vanishing}}{3.2.2}{}\ifthenelse{\equal{#1}{eqn--claim1}}{2.2.18}{}\ifthenelse{\equal{#1}{eqn--colim.DGCat.general}}{2.2.4}{}\ifthenelse{\equal{#1}{eqn--colim.DGCat}}{2.2.3}{}\ifthenelse{\equal{#1}{eqn--Det.calX}}{2.3.10}{}\ifthenelse{\equal{#1}{eqn--DM.calX}}{2.3.9}{}\ifthenelse{\equal{#1}{eqn--DM.Cech}}{2.2.20}{}\ifthenelse{\equal{#1}{eqn--DM.G.colim}}{2.3.5}{}\ifthenelse{\equal{#1}{eqn--DM.K-theory}}{2.1.9}{}\ifthenelse{\equal{#1}{eqn--DM.tau}}{2.2.17}{}\ifthenelse{\equal{#1}{eqn--DTM.asymmetric}}{5.3.3}{}\ifthenelse{\equal{#1}{eqn--Hom.X.Y}}{3.2.13}{}\ifthenelse{\equal{#1}{eqn--iota.Fl}}{5.0.2}{}\ifthenelse{\equal{#1}{eqn--iota.v.w}}{5.1.2}{}\ifthenelse{\equal{#1}{eqn--localization}}{2.1.5}{}\ifthenelse{\equal{#1}{eqn--localizationOne}}{2.1.4}{}\ifthenelse{\equal{#1}{eqn--Loop_Grp_Dfn}}{4.2.1}{}\ifthenelse{\equal{#1}{eqn--motivic.t.structure}}{3.2.5}{}\ifthenelse{\equal{#1}{eqn--presentation.Ind.Artin}}{2.3.1}{}\ifthenelse{\equal{#1}{eqn--presentation.Ind.scheme}}{2.4.1}{}\ifthenelse{\equal{#1}{eqn--relative.purity}}{2.1.8}{}\ifthenelse{\equal{#1}{eqn--rho.f!}}{2.3.8}{}\ifthenelse{\equal{#1}{eqn--VE}}{A.4.4}{}\ifthenelse{\equal{#1}{eqn--X.S.GX.etc}}{3.1.24}{}\ifthenelse{\equal{#1}{exam--affine.proj}}{2.2.13}{}\ifthenelse{\equal{#1}{exam--basic.WT}}{3.1.17}{}\ifthenelse{\equal{#1}{exam--BS.exam}}{3.2.3}{}\ifthenelse{\equal{#1}{exam--descent.Tate.not}}{3.1.25}{}\ifthenelse{\equal{#1}{exam--double.quot}}{4.2.12}{}\ifthenelse{\equal{#1}{exam--fusion.loop}}{6.1.3}{}\ifthenelse{\equal{#1}{exam--G.Whitney.Tate}}{3.1.13}{}\ifthenelse{\equal{#1}{exam--groups}}{A.4.12}{}\ifthenelse{\equal{#1}{exam--length.function.exam}}{4.2.16}{}\ifthenelse{\equal{#1}{exam--monoidal.unit}}{2.4.3}{}\ifthenelse{\equal{#1}{exam--motive.grass}}{2.3.13}{}\ifthenelse{\equal{#1}{exam--parahoric}}{4.2.2}{}\ifthenelse{\equal{#1}{exam--simple.reflection}}{4.3.14}{}\ifthenelse{\equal{#1}{exam--stratified.dfn}}{3.1.6}{}\ifthenelse{\equal{#1}{fiber_over_Grass}}{6.1.4}{}\ifthenelse{\equal{#1}{flag_act}}{4.3.2}{}\ifthenelse{\equal{#1}{foot.separated}}{1}{}\ifthenelse{\equal{#1}{funda.diag}}{6.2.3}{}\ifthenelse{\equal{#1}{glob_loop_group}}{6.1.1}{}\ifthenelse{\equal{#1}{intersection.complex}}{3.3.6}{}\ifthenelse{\equal{#1}{invariant.map}}{6.2.1}{}\ifthenelse{\equal{#1}{IW_Indentify}}{4.2.9}{}\ifthenelse{\equal{#1}{IW_Sub_Indentify}}{4.2.10}{}\ifthenelse{\equal{#1}{lemm--adm.finite.type}}{A.3.5}{}\ifthenelse{\equal{#1}{lemm--adm.reduced}}{A.3.3}{}\ifthenelse{\equal{#1}{lemm--adm.strata}}{A.3.2}{}\ifthenelse{\equal{#1}{lemm--affine.proj}}{4.2.7}{}\ifthenelse{\equal{#1}{lemm--base.change.schubert.field}}{4.3.6}{}\ifthenelse{\equal{#1}{lemm--compact.motives.Ind.Artin}}{2.3.6}{}\ifthenelse{\equal{#1}{lemm--descent.Beck.Chevalley}}{2.1.11}{}\ifthenelse{\equal{#1}{lemm--DM.double.tau}}{5.3.2}{}\ifthenelse{\equal{#1}{lemm--DM.G.BarC}}{2.2.7}{}\ifthenelse{\equal{#1}{lemm--DM.G/H}}{2.2.21}{}\ifthenelse{\equal{#1}{lemm--double.orbit}}{4.2.11}{}\ifthenelse{\equal{#1}{lemm--functoriality.equivariant}}{2.2.9}{}\ifthenelse{\equal{#1}{lemm--ind.pres}}{4.3.3}{}\ifthenelse{\equal{#1}{lemm--intermediate.simple}}{3.3.4}{}\ifthenelse{\equal{#1}{lemm--lim.equivalence}}{2.2.12}{}\ifthenelse{\equal{#1}{lemm--loc.triv}}{2.2.23}{}\ifthenelse{\equal{#1}{lemm--Lurie.co.limit}}{2.3.2}{}\ifthenelse{\equal{#1}{lemm--middle.extension}}{3.3.3}{}\ifthenelse{\equal{#1}{lemm--monadic.Ind}}{2.1.16}{}\ifthenelse{\equal{#1}{lemm--orbit.flag}}{4.3.7}{}\ifthenelse{\equal{#1}{lemm--parahoric.defi}}{4.2.4}{}\ifthenelse{\equal{#1}{lemm--pi*.fullyfaithful.MTM}}{3.2.12}{}\ifthenelse{\equal{#1}{lemm--pro.group}}{A.2.1}{}\ifthenelse{\equal{#1}{lemm--smooth.detects.t-structure}}{3.2.11}{}\ifthenelse{\equal{#1}{lemm--smooth.detects.Tate}}{3.1.20}{}\ifthenelse{\equal{#1}{lemm--t-structure.stratum}}{3.2.4}{}\ifthenelse{\equal{#1}{lemm--t.structure.limit}}{3.2.18}{}\ifthenelse{\equal{#1}{lemm--Tate.conservative}}{3.2.8}{}\ifthenelse{\equal{#1}{lemm--Tate.proper.descent}}{3.1.19}{}\ifthenelse{\equal{#1}{lemm--Tate.up.down}}{3.1.18}{}\ifthenelse{\equal{#1}{lemm--torsor.sequence}}{A.4.3}{}\ifthenelse{\equal{#1}{lemm--torsors}}{2.2.22}{}\ifthenelse{\equal{#1}{lemm--var.action}}{6.1.5}{}\ifthenelse{\equal{#1}{length_Grass_Iwahori}}{4.2.18}{}\ifthenelse{\equal{#1}{length_Grass}}{4.2.17}{}\ifthenelse{\equal{#1}{map.intro}}{1.2.1}{}\ifthenelse{\equal{#1}{middle-extension}}{3.3.2}{}\ifthenelse{\equal{#1}{nota--BS.vanishing.ladic}}{5.0.1}{}\ifthenelse{\equal{#1}{nota--BS.vanishing}}{3.2.1}{}\ifthenelse{\equal{#1}{nota--S.nochmal}}{3.0.1}{}\ifthenelse{\equal{#1}{nota--S}}{2.0.1}{}\ifthenelse{\equal{#1}{parahoric}}{4.2.3}{}\ifthenelse{\equal{#1}{prop--ballaballa}}{3.2.22}{}\ifthenelse{\equal{#1}{prop--boxtimes}}{2.4.4}{}\ifthenelse{\equal{#1}{prop--cells.flag}}{4.3.9}{}\ifthenelse{\equal{#1}{prop--change.facet}}{4.3.13}{}\ifthenelse{\equal{#1}{prop--DM.G.homotopy.invariant}}{2.2.11}{}\ifthenelse{\equal{#1}{prop--DTM.Fl.characterization}}{5.2.2}{}\ifthenelse{\equal{#1}{prop--DTM.G}}{3.1.27}{}\ifthenelse{\equal{#1}{prop--DTM.G/H}}{3.1.23}{}\ifthenelse{\equal{#1}{prop--equivariant.MTM}}{3.2.20}{}\ifthenelse{\equal{#1}{prop--existence.functors}}{6.3.3}{}\ifthenelse{\equal{#1}{prop--f_!.ind-Artin}}{2.3.3}{}\ifthenelse{\equal{#1}{prop--generators.DTM.G}}{3.2.23}{}\ifthenelse{\equal{#1}{prop--MTM.G}}{3.2.15}{}\ifthenelse{\equal{#1}{prop--Schubert.base.change}}{4.4.3}{}\ifthenelse{\equal{#1}{prop--sheafification.iso}}{2.2.25}{}\ifthenelse{\equal{#1}{prop--unipotent}}{A.4.6}{}\ifthenelse{\equal{#1}{prop--unseparated}}{2.1.14}{}\ifthenelse{\equal{#1}{prop--vector.extension}}{A.4.9}{}\ifthenelse{\equal{#1}{Quasi_Coxeter}}{4.2.13}{}\ifthenelse{\equal{#1}{reduced.eq}}{A.3.4}{}\ifthenelse{\equal{#1}{rema--bounded_subsets}}{4.2.6}{}\ifthenelse{\equal{#1}{rema--classical.equivariant}}{3.2.19}{}\ifthenelse{\equal{#1}{rema--classical.t-structure}}{3.2.7}{}\ifthenelse{\equal{#1}{rema--DM.prestacks}}{2.2.2}{}\ifthenelse{\equal{#1}{rema--DTM.large}}{3.1.9}{}\ifthenelse{\equal{#1}{rema--explain.stratified.Ind.scheme}}{3.1.3}{}\ifthenelse{\equal{#1}{rema--general.cellular}}{4.2.8}{}\ifthenelse{\equal{#1}{rema--l.adic.prestacks}}{2.3.11}{}\ifthenelse{\equal{#1}{rema--prestacks.examples}}{2.2.5}{}\ifthenelse{\equal{#1}{rema--realization.functor}}{3.2.9}{}\ifthenelse{\equal{#1}{rema--WT.maps}}{3.1.16}{}\ifthenelse{\equal{#1}{rema--WT.properties.fusion}}{6.3.7}{}\ifthenelse{\equal{#1}{rema--WT.smooth.duality}}{3.1.14}{}\ifthenelse{\equal{#1}{Root_Groups}}{4.3.10}{}\ifthenelse{\equal{#1}{Schubert_Map_Rel}}{4.4.2}{}\ifthenelse{\equal{#1}{Schubert_Map}}{4.3.5}{}\ifthenelse{\equal{#1}{Schubert_praesi_rel}}{4.4.4}{}\ifthenelse{\equal{#1}{Schubert_praesi}}{4.3.8}{}\ifthenelse{\equal{#1}{sect--algebraic.grps}}{A.2}{}\ifthenelse{\equal{#1}{sect--chevalley.base.change}}{4.4}{}\ifthenelse{\equal{#1}{sect--DM.Artin}}{2.3}{}\ifthenelse{\equal{#1}{sect--DM.ind-schemes}}{2.4}{}\ifthenelse{\equal{#1}{sect--DM.prestacks}}{2.2}{}\ifthenelse{\equal{#1}{sect--DM.schemes}}{2.1}{}\ifthenelse{\equal{#1}{sect--DM}}{2}{}\ifthenelse{\equal{#1}{sect--DTM.definitions}}{3.1}{}\ifthenelse{\equal{#1}{sect--DTM.double}}{5.3}{}\ifthenelse{\equal{#1}{sect--DTM.Fl}}{5}{}\ifthenelse{\equal{#1}{sect--DTM}}{3}{}\ifthenelse{\equal{#1}{sect--ind.schemes}}{A.1}{}\ifthenelse{\equal{#1}{sect--intersection}}{6}{}\ifthenelse{\equal{#1}{sect--invariant}}{6.2}{}\ifthenelse{\equal{#1}{sect--loop.definitions}}{4.2}{}\ifthenelse{\equal{#1}{sect--loop.group.dfn}}{4.1}{}\ifthenelse{\equal{#1}{sect--loop.grps}}{4}{}\ifthenelse{\equal{#1}{sect--MTM}}{3.2}{}\ifthenelse{\equal{#1}{sect--pro.action}}{A.3}{}\ifthenelse{\equal{#1}{sect--realization}}{2.1.2}{}\ifthenelse{\equal{#1}{sect--shtukas.def}}{6.3}{}\ifthenelse{\equal{#1}{sect--Stratifications.flag}}{4.3}{}\ifthenelse{\equal{#1}{sect--Tate.Fl}}{5.2}{}\ifthenelse{\equal{#1}{sect--torsors}}{A.4}{}\ifthenelse{\equal{#1}{sect--WT.Fl}}{5.1}{}\ifthenelse{\equal{#1}{sheaf_condition}}{2.2.14}{}\ifthenelse{\equal{#1}{shtuka.inv}}{6.3.2}{}\ifthenelse{\equal{#1}{shtuka.pres}}{6.3.1}{}\ifthenelse{\equal{#1}{strat.map.square}}{3.1.2}{}\ifthenelse{\equal{#1}{stratified.diagram}}{3.1.7}{}\ifthenelse{\equal{#1}{stratified.Ind-pres}}{3.1.4}{}\ifthenelse{\equal{#1}{syno--motives}}{2.1.1}{}\ifthenelse{\equal{#1}{Tate.big.orbit}}{5.1.3}{}\ifthenelse{\equal{#1}{theo--D.et!.sheaf}}{2.1.15}{}\ifthenelse{\equal{#1}{theo--descent.prestacks}}{2.2.16}{}\ifthenelse{\equal{#1}{theo--DM.descent}}{2.1.13}{}\ifthenelse{\equal{#1}{theo--equivariant.Chow}}{2.2.10}{}\ifthenelse{\equal{#1}{theo--equivariant.DTM.flag}}{5.3.4}{}\ifthenelse{\equal{#1}{theo--f!.Artin}}{2.3.7}{}\ifthenelse{\equal{#1}{theo--Fl.WT}}{5.1.1}{}\ifthenelse{\equal{#1}{theo--generators.DTM.flag}}{5.2.3}{}\ifthenelse{\equal{#1}{theo--motives.Ind-schemes}}{2.4.2}{}\ifthenelse{\equal{#1}{theo--simple.objects}}{3.3.8}{}\ifthenelse{\equal{#1}{tocindent-1}}{0pt}{}\ifthenelse{\equal{#1}{tocindent0}}{0pt}{}\ifthenelse{\equal{#1}{tocindent1}}{66.11127pt}{}\ifthenelse{\equal{#1}{tocindent2}}{0pt}{}\ifthenelse{\equal{#1}{tocindent3}}{0pt}{}\ifthenelse{\equal{#1}{torsor.lim}}{A.4.1}{}\ifthenelse{\equal{#1}{torsor.lim1}}{A.4.2}{}\ifthenelse{\equal{#1}{unif.unab}}{6.1.6}{}\ifthenelse{\equal{#1}{vector_bundle}}{A.4.11}{}\ifthenelse{\equal{#1}{Weil_resitriction}}{4.2.5}{}}
\def \inter#1{\cite[#1]{RicharzScholbach:Intersection}}
\def \irefco#1{\inter{Corollary~\iref{coro--#1}}}
\def \irefde#1{\inter{Definition~\iref{defi--#1}}}
\def \irefdele#1{\inter{Definition and Lemma~\iref{defilemm--#1}}}

\def \irefle#1{\inter{Lemma~\iref{lemm--#1}}}
\def \irefpr#1{\inter{Proposition~\iref{prop--#1}}}
\def \irefre#1{\inter{Remark~\iref{rema--#1}}}
\def \irefth#1{\inter{Theorem~\iref{theo--#1}}}
\def \irefsect#1{\inter{§\iref{sect--#1}}}
\def \irefsy#1{\inter{Synopsis~\iref{syno--#1}}}

\usepackage{ifthen}

\def\satref#1{\ifthenelse{\equal{#1}{associative}}{Lemma~3.7}{}\ifthenelse{\equal{#1}{base.change.nice}}{Lemma~2.12}{}\ifthenelse{\equal{#1}{base.change.t.exact}}{Proposition~2.15}{}\ifthenelse{\equal{#1}{base.scheme}}{Notation~2.1}{}\ifthenelse{\equal{#1}{constraints}}{Proposition~5.9}{}\ifthenelse{\equal{#1}{convolution.DM.ell}}{Proposition~3.14}{}\ifthenelse{\equal{#1}{convolution.DTM.weights}}{Proposition~4.8}{}\ifthenelse{\equal{#1}{convolution.Fl.Tate.2x.Iwahori}}{Proposition~3.26}{}\ifthenelse{\equal{#1}{convolution.Fl.Tate.Iwahori}}{Proposition~3.19}{}\ifthenelse{\equal{#1}{convolution.Fl.Tate}}{Theorem~3.17}{}\ifthenelse{\equal{#1}{convolution.HoDM.independent}}{Proposition~3.4}{}\ifthenelse{\equal{#1}{convolution.product}}{Definition~3.1}{}\ifthenelse{\equal{#1}{convolution.Tate.perverse}}{Lemma~5.8}{}\ifthenelse{\equal{#1}{Day.convolution}}{Lemma~A.2}{}\ifthenelse{\equal{#1}{DM!.placid.prestacks}}{Corollary~A.9}{}\ifthenelse{\equal{#1}{DM*.PreStk.monoidal}}{Corollary~A.3}{}\ifthenelse{\equal{#1}{DM*.slm}}{Lemma~A.1}{}\ifthenelse{\equal{#1}{DM*vs!.prestacks}}{Corollary~A.8}{}\ifthenelse{\equal{#1}{DM*vs!}}{Proposition~A.4}{}\ifthenelse{\equal{#1}{epsilon.remark}}{Remark~6.7}{}\ifthenelse{\equal{#1}{equivalence}}{Corollary~5.7}{}\ifthenelse{\equal{#1}{explain.convolution}}{Remark~3.2}{}\ifthenelse{\equal{#1}{extensions.grass}}{Corollary~5.5}{}\ifthenelse{\equal{#1}{exterior.product.weights}}{Lemma~4.3}{}\ifthenelse{\equal{#1}{f!.nice}}{Theorem~2.14}{}\ifthenelse{\equal{#1}{f*.convolution}}{Lemma~3.15}{}\ifthenelse{\equal{#1}{fiber.functor}}{Definition~5.11}{}\ifthenelse{\equal{#1}{finite.field.exam}}{Example~6.12}{}\ifthenelse{\equal{#1}{full.Tannaka}}{Theorem~6.14}{}\ifthenelse{\equal{#1}{Hom.IC}}{Lemma~6.3}{}\ifthenelse{\equal{#1}{IC.pure}}{Theorem~4.1}{}\ifthenelse{\equal{#1}{IC.star.skyscraper}}{Lemma~5.4}{}\ifthenelse{\equal{#1}{KL.vanishing.theo}}{Theorem~2.21}{}\ifthenelse{\equal{#1}{motivation.formula}}{Remark~3.24}{}\ifthenelse{\equal{#1}{motivic.t.structure}}{Theorem~2.8}{}\ifthenelse{\equal{#1}{MTM.basics}}{Synopsis~5.1}{}\ifthenelse{\equal{#1}{orbit.simply.conn}}{Lemma~2.9}{}\ifthenelse{\equal{#1}{parity.coro.grass.last}}{Corollary~5.12}{}\ifthenelse{\equal{#1}{parity.coro.grass}}{Corollary~5.3}{}\ifthenelse{\equal{#1}{parity.decomposition}}{Corollary~5.6}{}\ifthenelse{\equal{#1}{pi*.equivalence.MTM}}{Lemma~2.10}{}\ifthenelse{\equal{#1}{placid}}{Definition~A.5}{}\ifthenelse{\equal{#1}{prestack.DM}}{Definition~2.4}{}\ifthenelse{\equal{#1}{properties.DM}}{Theorem~2.6}{}\ifthenelse{\equal{#1}{pushforward.symmetric.monoidal}}{Proposition~5.13}{}\ifthenelse{\equal{#1}{realization.functor.coro}}{Corollary~2.20}{}\ifthenelse{\equal{#1}{realization.functor}}{Theorem~2.19}{}\ifthenelse{\equal{#1}{Satake.convolution}}{Lemma~6.5}{}\ifthenelse{\equal{#1}{Satake.Ind.objects}}{Corollary~6.15}{}\ifthenelse{\equal{#1}{Satake}}{Theorem~6.8}{}\ifthenelse{\equal{#1}{star.twiddle}}{Lemma~3.11}{}\ifthenelse{\equal{#1}{Tannaka_Cat}}{Theorem~5.14}{}\ifthenelse{\equal{#1}{Tannaka.Sat}}{Corollary~6.6}{}\ifthenelse{\equal{#1}{twiddle.basic}}{Lemma~3.12}{}\ifthenelse{\equal{#1}{weight.creation}}{Lemma~2.17}{}\ifthenelse{\equal{#1}{weight.grading}}{Corollary~4.2}{}\ifthenelse{\equal{#1}{weights.DM.G}}{Definition~4.5}{}\ifthenelse{\equal{#1}{classical.satake.iso}}{1.1}{}\ifthenelse{\equal{#1}{compatibility.realization.convolution}}{3.1.2}{}\ifthenelse{\equal{#1}{convolution.Fl.Tate:3}}{3.23}{}\ifthenelse{\equal{#1}{DM.intro}}{1.2}{}\ifthenelse{\equal{#1}{eqn--base.change.DTM}}{2.11}{}\ifthenelse{\equal{#1}{eqn--boxy}}{3.3}{}\ifthenelse{\equal{#1}{eqn--convolution.Fl.Tate:2}}{3.21}{}\ifthenelse{\equal{#1}{eqn--convolution.terms}}{3.13}{}\ifthenelse{\equal{#1}{eqn--convolution}}{3.10}{}\ifthenelse{\equal{#1}{eqn--define.convolution}}{1.3}{}\ifthenelse{\equal{#1}{eqn--DM.colimitissimo}}{4.7}{}\ifthenelse{\equal{#1}{eqn--DM.prestacks}}{2.5}{}\ifthenelse{\equal{#1}{eqn--DM.Sch.ft}}{2.3}{}\ifthenelse{\equal{#1}{eqn--f*.IC}}{2.16}{}\ifthenelse{\equal{#1}{eqn--iso.Fls.proof}}{3.20}{}\ifthenelse{\equal{#1}{eqn--iso.Flw.proof}}{3.22}{}\ifthenelse{\equal{#1}{eqn--Iwahori}}{3.16}{}\ifthenelse{\equal{#1}{eqn--MTM.rl}}{3.5}{}\ifthenelse{\equal{#1}{eqn--Sat.MTM}}{6.1}{}\ifthenelse{\equal{#1}{eqn--star.associative}}{3.8}{}\ifthenelse{\equal{#1}{eqn--X.lim}}{A.6}{}\ifthenelse{\equal{#1}{full.tannaka.group}}{6.13}{}\ifthenelse{\equal{#1}{functions.sec}}{6.4}{}\ifthenelse{\equal{#1}{intersection.intro}}{1.4}{}\ifthenelse{\equal{#1}{iota.vw}}{2.13}{}\ifthenelse{\equal{#1}{iso.groups.tori}}{6.10}{}\ifthenelse{\equal{#1}{iso.groups}}{6.9}{}\ifthenelse{\equal{#1}{MTM.general.intro}}{1.5}{}\ifthenelse{\equal{#1}{orbit.inclusion}}{2.2}{}\ifthenelse{\equal{#1}{parity.eqn}}{2.22}{}\ifthenelse{\equal{#1}{parity.fct}}{5.2}{}\ifthenelse{\equal{#1}{really.to.show}}{3.18}{}\ifthenelse{\equal{#1}{reformulation.convolution}}{3.1.1}{}\ifthenelse{\equal{#1}{sect--box.product}}{A}{}\ifthenelse{\equal{#1}{sect--changing.base}}{2.4}{}\ifthenelse{\equal{#1}{sect--convolution.product}}{3}{}\ifthenelse{\equal{#1}{sect--defi.associ}}{3.1}{}\ifthenelse{\equal{#1}{sect--DTM.Tate}}{3.2}{}\ifthenelse{\equal{#1}{sect--dual.grp.sect}}{6}{}\ifthenelse{\equal{#1}{sect--loop.grps.Satake}}{2.1}{}\ifthenelse{\equal{#1}{sect--motive.affine.flag.basics}}{2}{}\ifthenelse{\equal{#1}{sect--motives.prestack}}{2.2}{}\ifthenelse{\equal{#1}{sect--parity.vanishing}}{2.5}{}\ifthenelse{\equal{#1}{sect--purity}}{4}{}\ifthenelse{\equal{#1}{sect--Satake.category}}{6.1}{}\ifthenelse{\equal{#1}{sect--stratified.motives}}{2.3}{}\ifthenelse{\equal{#1}{sect--tate.aff.grass}}{5}{}\ifthenelse{\equal{#1}{sect--tensor.structure.mtm}}{5.2}{}\ifthenelse{\equal{#1}{trace.eq}}{6.16}{}}

\def\sref#1{\cite[\satref{#1}]{RicharzScholbach:Motivic}}

\begin{document}

\author{Timo Richarz and Jakob Scholbach}
\title{Tate motives on Witt vector affine flag varieties}

\begin{abstract}
Relying on recent advances in the theory of motives we develop a general formalism for derived categories of motives with $\bbQ$-coefficients on perfect \ii-prestacks.
We construct Grothendieck's six functors for motives over perfect (ind-)schemes perfectly of finite presentation.  
Following ideas of Soergel--Wendt, this is used to study basic properties of stratified Tate motives on Witt vector partial affine flag varieties.
As an application we give a motivic refinement of Zhu's geometric Satake equivalence for Witt vector affine Grassmannians in this set-up.
\end{abstract}


\maketitle

\tableofcontents

\section{Introduction}

\subsection{Motivation}
\label{sect--motivation}
The Satake isomorphism (\cite{Gross:Satake}) plays a foundational rôle in several branches of the Langlands program.
Let $p$ be a prime number, and fix $p^{1/2}\in \overline\bbQ$.
For a split reductive group $G$ and a non-archimedean local field like $\bbQ_p$ or $\bbF_p\rpot{t}$, this isomorphism relates spherical functions pertaining to $G$ to the representation theory of the Langlands dual group $\widehat G$:
$$\calC^\comp\big(G(\Qp) /\!/ G(\Zp)\big)\; \stackrel \cong\lr\; R(\widehat G) \;\stackrel \cong \longleftarrow\; \calC^\comp\big(G(\Fp\rpot t) /\!/ G(\Fp \pot t)\big)\eqlabel{intro.iso}$$
The so-called spherical Hecke algebras at the right and left consist of finitely supported $\bbZ[p^{\pm 1/2}]$-valued functions on the $G(\Zp)$- (resp.~$G(\Fp \pot t)$-)double cosets in $G(\Qp)$ (resp. $G(\Fp \rpot t)$).
In the middle, $\widehat G$ denotes the Langlands dual group formed over $\bbQ$ using a fixed pinning $(T, B, X)$ of $G$.
By the classification of split reductive groups in terms of their root data, it can be described by switching characters and cocharacters, as well as roots and coroots.  
For example (\cite{Cogdell:DualGroup}), $\widehat{\GL}_n=\GL_n$, $\widehat{\SL}_n = \PGL_n$ and $\widehat{\on{Sp}}_{2g}=\on{SO}_{2g+1}$.
The ring $R(\widehat G)$ is the Grothendieck $\bbZ[p^{\pm 1/2}]$-algebra of algebraic representations of this dual group.

The chain of isomorphisms \refeq{intro.iso} is an instance of the resemblance (\cite{Weil:Analogie}) between $\Qp$ and $\Fp \rpot t$, in which the uniformizer $p$ corresponds to the uniformizer $t$.
Despite being only a superficial similarity --given that the addition and multiplication in the two fields are severely different-- the kinship is strong enough so that the characteristic functions of the double cosets $G(\Zp)p^\mu G(\Zp)$ and $G(\Fp \pot t)t^\mu G(\Fp \pot t)$ correspond to each other under \refeq{intro.iso}.
Here $\mu\co \Gm \r T$ is any dominant cocharacter.

The geometric Satake equivalence might be regarded, at the same time, as a geometrization and a categorification of the right hand isomorphism above: the right-most term is geometrized by means of the affine Grassmannian, whose $\bbF_p$-points are the quotient $G(\Fp \rpot t) / G(\Fp \pot t)$. The categorification is achieved by relating the entire category of algebraic $\widehat G$-representations (as opposed to its Grothendieck ring) to certain equivariant perverse sheaves on the affine Grassmannian.
This circle of ideas has been worked on by various authors including Lusztig \cite{Lusztig:Singularities}, Ginzburg \cite{Ginzburg:Perverse}, Belinson--Drinfeld \cite{BeilinsonDrinfeld:Quantization} and Mirkovi\'c--Vilonen \cite{MirkovicVilonen:Geometric}. 
For the relation with \refeq{intro.iso} the reader is referred to \cite{RicharzZhu:Ramified, Zhu:Introduction}.

In analogy to the function field case, a geometrization and categorification of the left hand isomorphism in \refeq{intro.iso} was proven by Zhu \cite{Zhu:Affine}. 
To begin with, such an endeavour is made possible by \cite{Zhu:Affine} and Bhatt--Scholze \cite{BhattScholze:Projectivity} who have shown that the cosets $G(\Qp) / G(\Zp)$ are the $\bbF_p$-points of an algebro-geometric object known as the Witt vector affine Grassmannian, and denoted $\Gr_G\to \Spec(\bbF_p)$ in this paper.
For a fixed prime nubmer $\ell\not = p$ and a choice of half twist $\overline{\bbQ}_\ell({1\over 2})$, Zhu proves an equivalence of Tannakian categories
$$\Perv_{L^{+} G}\big(\Gr_{G}; \overline{\bbQ}_\ell\big)^{0,\on{ss}}\;\stackrel \cong \lr \; \Rep_{\overline{\bbQ}_\ell}\big(\widehat G\big) \eqlabel{intro.equivalence}$$
between semi-simple perverse $\overline \bbQ _\ell$-sheaves on $\Gr_{G}$ of weight $0$ which are furthermore equivariant with respect to the action of the positive loop group $L^+G$.
This positive loop group has as $\Fp$-points the group $G(\Zp)$; it lies inside the full loop group $LG$, whose $\Fp$-points are $G(\Qp)$.

Building upon our earlier work \cite{RicharzScholbach:Intersection, RicharzScholbach:Motivic}, the goal of the present paper is to refine \refeq{intro.equivalence} into an equivalence between the category of mixed Tate motives on the double quotient $L^+G\bsl LG/ L^+G$ and representations of Deligne's extended dual group $\widehat G_1$, a certain extension of $\Gm$ (which records the weights) by $\widehat G$. 
This equivalence is independent of the choice of $\ell\not = p $ present in \refeq{intro.equivalence} via the use of $\ell$-adic cohomology.
We refer to \myref{classical.satake.rema} for the relation with \refeq{intro.iso}.
Also we remark that Zhu \cite{Zhu:Geometric} has explained the construction of a motivic Satake equivalence using numerical motives.
This approach is based on an explicit enumeration of algebraic cycles on affine Grassmannians. 
By comparison, the approach taken in this paper is more strongly relying on the general framework of motives given by Ayoub \cite{Ayoub:Six1, Ayoub:Six2} and Cisinski--D\'{e}glise \cite{CisinskiDeglise:Triangulated}.

\subsection{Results}
A prestack over $\bbF_p$ is a functor from affine $\bbF_p$-algebras to anima (called spaces in \cite{Lurie:Higher}, Kan complexes, or \ii-groupoids or, in classical terminology, simplicial sets up to weak equivalences).
Examples of prestacks include (ordinary) presheaves such as $\bbF_p$-schemes and ind-schemes (via their functor of points), but also homotopy quotients of such objects by group actions, and more general geometric objects such as higher stacks over $\bbF_p$.
In \cite{RicharzScholbach:Intersection} we have constructed a stable \ii-category $\DM(X)$ of motives with $\bbQ$-coefficients on any prestack $X$.
For example, we can speak of and conveniently work with $\DM(L^+G \bsl LG / L^+G)$, the category of motives on the indicated double quotient.
 
In order to state our results, we first need to deal with the general process of {\it perfection}, see \S6 in the first arXiv version of \cite{RicharzScholbach:Intersection}, and also \cite{ElmantoKhan:Perfection}, or the recent preprint of Bouthier--Kazhdan--Varshavsky \cite[\S2.3.5]{BouthierKazhdanVarshavsky:Perverse}.
There are two ways to turn a $\bbF_p$-prestack $X$ perfect, i.e., to ensure that pullback along the absolute Frobenius $\sigma$ is an isomorphism: the limit perfection
$$\on{lim}_\sigma X \defined \lim \big(\dots \overset{\sigma}{\r} X\overset{\sigma}{\r} X\overset{\sigma}{\r} X\big),$$
and the colimit perfection $X^\perf:=\on{colim}_{T\to X} \lim_\sigma T$.
The latter is constructed by glueing the perfections of all affine $\bbF_p$-schemes covering the prestack (see \myref{colimit.perfection}).

\prop \textup{(}\myref{tomaten.auf.den.augen}\textup{)}
\mylabel{intro.prop}
For any $\bbF_p$-prestack $X$, there are maps
$$\on{lim}_\sigma X \longleftarrow (\on{lim}_\sigma X)^\perf \stackrel \simeq \lr X^\perf$$
inducing equivalences on categories of motives
$$\DM(\on{lim}_\sigma X) \;\simeq\; \DM((\on{lim}_\sigma X)^\perf) \;\simeq\; \DM(X^\perf).$$
\xprop

Recall from \cite{Zhu:Affine,BhattScholze:Projectivity} (or see \refsect{motives.perfect.schemes} for a brief recapitulation) that perfectly finitely presented (pfp) $\bbF_p$-schemes are precisely those $\bbF_p$-schemes $X$ which arise as $X \simeq \lim_\sigma X_0$ for some finite-type $\bbF_p$-scheme $X_0$.
Such a scheme $X_0$ is called a model of $X$. We have $\DM(X_0) \simeq \DM(X)$ by the previous proposition. 
A typical example of a pfp scheme is $(\A^1_{\bbF_p})^\perf = \Spec \Fp[t^{p^{-\infty}}]$, the perfection of $\A^1_{\bbF_p}$.

\theo \textup(\myref{DM.perfect.schemes,functor.adic}\textup)
\mylabel{intro.theorem1}
The equivalences of \myref{intro.prop} can be used to endow motives on pfp schemes over $\bbF_p$ with a six functor formalism obeying all the standard properties including base change, localization, h-descent, duality and homotopy invariance with respect $(\bbA^1_{\bbF_p})^\perf$.
\xtheo  

Similar results also hold for motives on ind-(pfp schemes), i.e., objects presented as $\N$-indexed colimits $X = \colim X_i$, where the $X_i$ are pfp schemes and the transition maps $X_i \r X_{i+1}$ are closed immersions.
An example is the Witt vector affine Grassmannian $\Gr_G=(LG/L^+G)^\et$ which is shown to be an ind-(pfp scheme) in \cite{BhattScholze:Projectivity}.

Such a convenient formalism lends itself to applications in geometric representation theory. 
Extending previous work of Soergel--Wendt \cite{SoergelWendt:Perverse} for finite type schemes and of \cite{RicharzScholbach:Intersection} for ind-schemes of ind-(finite type) to the (ind-)pfp case, we say that a stratification on an ind-(pfp scheme) $X$ is Whitney--Tate if the derived categories of Tate motives on the individual strata ``glue'' in a meaningful manner (see \myref{Whitney-Tate.pfp} for the precise condition). 

\theo\textup(special case of \myref{Fl.MTM.adic}\textup)
\mylabel{intro.theorem2}
\begin{enumerate}
\item The stratification on the Witt vector affine Grassmannian $\Gr_G$ by $L^+G$-orbits is Whitney--Tate.
Accordingly, the stable \ii-category $\DTM(\Gr_G)$ of stratified Tate motives is well-defined: such motives are characterized by the property that their restriction to the $L^+G$-orbits are Tate, i.e., generated (under shifts and extensions) by motives of the form $1(n)$, $n\in \bbZ$ on the $L^+G$-orbits.

\item The category $\DTM(\Gr_G)$ admits a ``motivic'' $t$-structure whose heart $\MTM(\Gr_G)$ is generated by \emph{intersection motives} $\IC_\mu(n)$, $n\in \bbZ$ as simple objects, where $\mu\co \Gm \to T$ ranges over the dominant cocharacters of the fixed maximal torus $T\subset G$.
\end{enumerate}
\xtheo

In fact, this result is as a special case of \myref{Fl.MTM.adic} where we prove the same result more generally for stratifications on Witt vector (partial) affine flag varieties $\Fl_\bbf=(LG/\calP_\bbf)^\et$ associated to any parahoric subgroup $\calP_\bbf \subset L G$, in place of $L^+G$.
An example is the Iwahori subgroup $\calB\subset LG$ associated to the choice of the Borel subgroup $B\subset G$.

The convolution product of functions in the Hecke algebra is enhanced by the convolution product (\refsect{convolution.product.witt})
$$\str\star\str \co \DM(L^+G \bsl LG / L^+ G) \x \DM(L^+G \bsl LG / L^+ G) \r \DM(L^+G \bsl LG / L^+ G).$$
By étale descent, motives on $LG / L^+ G$ are equivalent to motives on the étale sheafification $\Gr_G = (LG / L^+ G)^\et$.
Therefore, motives on the double quotient $L^+G \bsl LG / L^+ G$ are equivalent to $L^+ G$-equivariant motives on $\Gr_G$.

We define a category of equivariant stratified Tate motives 
$$\DTM(L^+ G \backslash LG / L^+G) \subset \DM(L^+ G \bsl LG / L^+ G) \eqlabel{intro.dtm}$$ 
by requiring that the underlying (non-equivariant) motive on $\Gr_G$ is stratified Tate in the sense above.
By \myref{DTM.doubly.equivariant}, the category $\DTM(L^+ G \backslash LG / L^+G)$ inherits a t-structure from the t-structure on $\DTM(\Gr_G)$.
The heart of this t-structure is denoted $\MTM(L^+ G \backslash LG / L^+G)$.
In relation to the outline in \refsect{motivation}, this category is an independent-of-$\ell$ version of perverse $L^+G$-equivariant $\ell$-adic sheaves on $\Gr_G$ 
subject to the condition that their restriction to the strata are Tate twists of the (possibly shifted) constant sheaf.

\theo 
\textup(special case of \myref{convolution.product.witt}\textup)
The convolution product functor $\star$ preserves the subcategory \textup{\refeq{intro.dtm}}. 
\xtheo

Again, this statement is shown in \myref{convolution.product.witt} in the greater generality for Witt vector (partial) affine flag varieties as above.
The following statements, however, are specific to the choice of $L^+ G$.
In \myref{semisimplicity.prop} we prove, as a consequence of the semi-simplicity of $\MTM(\Fq)$ and Kazhdan--Lusztig's parity vanishing, that the forgetful functor 
$$\MTM(L^+ G \bsl LG / L^+ G) \r \MTM(\Gr_G) \eqlabel{intro.mtm}$$
is an equivalence of semi-simple (abelian) categories whose generators are the simple objects $\IC_\mu(n)$ mentioned in \myref{intro.theorem2} ii).
This category is stable under the convolution product $\star$ and hence inherits a monoidal structure by \myref{convolution.product.witt} iii).

We write $\widehat G_1$ for the modification of the Langlands dual group introduced by Frenkel and Gross \cite{FrenkelGross:Rigid} following a suggestion of Deligne \cite{Deligne:Letter2007}.
If $G$ is simply connected, $\widehat G_1$ is the (direct) product of $\widehat G$ and $\Gm$, where the $\Gm$-factor arises from the presence of Tate twists.
In general, $\widehat G_1$ is a not necessarily split extension of $\Gm$ by $\widehat G$, a phenomenon which is related to the necessity of fixing $p^{1/2}\in \overline\bbQ$ in the Satake isomorphisms \refeq{intro.iso}.
Our main result is the motivic Satake equivalence for Witt vector affine flag varieties. 

\theo \textup(\myref{motivic.satake.witt,Satake.equiv.thm}\textup)
The monoidal structure on \textup{\refeq{intro.mtm}} given by the convolution product is part of the structure of a $\bbQ$-linear Tannakian category, and as such there is, after extending scalars to $\overline \bbQ$, an equivalence
$$\MTM (L^+G \bsl LG / L^+ G;\overline\bbQ)\;\stackrel \simeq \lr\; \Rep_{\overline \bbQ}(\widehat G_1). $$
\xtheo

This statement is independent of $\ell$ since no category of $\ell$-adic (perverse) sheaves is required for its formulation.
Furthermore, this equivalence is uniquely determined by its compatibility with (the appropriate version of) the $\ell$-adic equivalence \refeq{intro.equivalence}.
As mentioned in \refsect{motivation} a similar result was also obtained in \cite{Zhu:Geometric} using a different route, namely an explicit enumeration of algebraic cycles on the Witt vector affine Grassmannian, and working with the category of numerical motives.
As is explained by the material in the sections below, working with the six functor formalism of motives as in this paper offers more flexibility including, say, handling the derived categories $\DTM(\calB \bsl LG / \calB)$ of doubly Iwahori-equivariant stratified Tate motives on the loop group: the convolution product on this category does not respect the abelian subcategory of equivariant mixed Tate motives, so that working with numerical motives is prohibitive in such a situation.

\bigskip
\noindent \textit{Acknowledgements.} 
We thank Eugen Hellmann and Thomas Nikolaus for helpful discussions, and the anonymous referee for many suggestions that improved the quality of the manuscript.
The authors thank the University of M\"unster, the Deutsche Forschungsgemeinschaft
(DFG, German Research Foundation; under Germany's Excellence Strategy EXC 2044–390685587, Mathematics Münster: Dynamics–Geometry–Structure) and the Technical University of Darmstadt for financial and logistical support which made this research possible.

\section{Motives on perfect (ind-)schemes}

\subsection{Recollections on prestacks}
\label{sect--DM.prestacks.recollections}
In order to treat motives on the Witt vector affine Grassmannian and to relate them to motives on ordinary affine Grassmannians, it is very useful to use the definition of a category of motives on arbitrary prestacks developed in \irefsect{DM}. 
We briefly recall the rudiments in this subsection.
A very similar theory, for $\ell$-adic sheaves instead of motives, has been developed very recently in \cite{BouthierKazhdanVarshavsky:Perverse}.

For a field $k$, $\AffSch_k^\ft$ is the category of affine $k$-schemes of finite type.
We denote by $\AffSch_k$ its $\kappa$-pro-completion, i.e., the category of affine $k$-schemes whose underlying $k$-algebra is generated by at most $\kappa$ elements. 
Here $\kappa$ is a large enough regular cardinal fixed once and for all.
For all purposes in this paper, we may choose $\kappa$ to be the countable cardinal.
The \ii-category of \emph{prestacks} is defined as the \ii-categorical presheaf category, i.e., 
$$\PreStk_k \defined \Fun\big(\AffSch_k^\opp, \Ani\big).$$
Here, following the terminology and discussion in \cite[§5.1.4]{CesnaviciusScholze:Purity}, $\Ani$ denotes the \ii-category of anima, i.e., the free completion of the category of finite sets under sifted homotopy colimits, or, equivalently, the \ii-category of spaces \cite[§1.2.16]{Lurie:Higher}, also called \ii-groupoids.
Prestacks are an extremely general class of algebro-geometric objects: by means of its functor of points, any $k$-scheme defines a prestack.
In addition, prestacks are closed under all (homotopy) limits and colimits.
In particular, any strict ind-scheme 
$X = \colim X_i$ is a prestack, as is a quotient (always understood in the \ii-categorical sense, otherwise also known as a homotopy quotient) $Y / G$, where $G$ is a group-valued presheaf acting on a presheaf $Y$.

We denote by $\DGCat_\cont$ the $\infty$-category of presentable, stable dg-$\infty$-categories with continuous (i.e., colimit-preserving) functors.

\cons
\mylabel{DM.prestacks.construction}
In \irefde{DM.prestacks}, we have defined a category $\DM(X)$ in $\DGCat_\cont$ of \emph{motives on any prestack $X$} with the following features:
\begin{enumerate}
\item
For a scheme $X$ of finite type over $k$, the homotopy category of $\DM(X)$ (which is a triangulated category) recovers the usual triangulated category of motives with rational coefficients due to Ayoub and Cisinski--Déglise \cite{Ayoub:Six1,Ayoub:Six2}, \cite[§14]{CisinskiDeglise:Triangulated}.

\item
If $X = \lim_{n \in \N} X_n$ is an affine pro-algebraic scheme, i.e., $X_i / k$ is affine and of finite type and so are the transition maps $X_n \stackrel{p_n} \lr X_{n-1}$, then 
$$\DM(X) = \colim \left ( \DM(X_0) \stackrel{(p_0)^!} \lr \DM(X_1) \stackrel{(p_1)^!} \lr  \dots \right )$$
where the colimit is taken in $\DGCat_\cont$ (we emphasize that this category consists of colimit-preserving functors only).

\item
If a prestack $X$ admits a presentation $X = \colim X_i$, where $X_i \in \AffSch_k$, then
$$\DM(X) = \lim \big(\DM(X_i), (f_{ij})^!\big)$$
where the maps $X_i \stackrel{f_{ij}}\lr X_j$ are the ones in the colimit diagram defining $X$.
More colloquially speaking, a motive on $X$ is a family of motives on all the ``covers'' $X_i$ compatible with the system defining the colimit diagram.

\item
\label{item--DM!.construction}
For any map of prestacks $f\co X \r Y$, there is a functor $f^!\co \DM(Y) \r \DM(X)$. 
These yield a presheaf $\DM\co \PreStk^\opp \r \DGCat_\cont$ which is shown to be a sheaf in Voevodsky's $h$-topology in \irefth{descent.prestacks}.
If we need to emphasize the nature of the pullback functors, we also write $\DM^!$ for this (pre)sheaf.

\item
In both preceding items, we use $!$-pullback instead of the more easily defined $*$-pullback in order to obtain the desired categories $\DM(X)$ for ind-schemes $X$, see \cite[Corollary~2.3.4]{RicharzScholbach:Intersection}.
For $X = Y / G$, the quotient of an action of a smooth algebraic group $G$ on a finite type $k$-scheme $Y$, the approach with $*$- and $!$-pullback yields equivalent categories, though, cf.~the discussion in \irefre{DM.prestacks} and after \irefle{DM.G.BarC}.
\end{enumerate}
\xcons

\subsection{Perfection of prestacks}
In this section we develop a general formalism of colimit perfections for \ii-prestacks, and compare it to the limit perfection as used in \cite[Appendix~A]{Zhu:Affine}, \cite{BhattScholze:Projectivity} (cf.~also \cite{BertapelleGonzalez:Perfection}) and \cite[Appendix~A]{XiaoZhu:Cycles}.

Let $p$ be a prime, and let $k$ be a perfect field of characteristic $p$. For each $T = \Spec R\in \AffSch_k$, let $\sigma_T\co T\to T$ be the absolute Frobenius morphism given by $R\to R, x\mapsto x^p$. This defines an endofunctor $\sigma$ of the identity transformation on $\AffSch_k$. Let $\AffSch_k^{\perf}$ be the full subcategory of objects $T\in\AffSch_k$ such that $\sigma_T$ is an automorphism. These are called {\em perfect affine $k$-schemes}. The inclusion $\AffSch_k^\perf\subset \AffSch_k$ admits a right adjoint given by $T\mapsto \lim_\sigma T=\Spec(\colim_{r\mapsto r^p}R)$.

Precomposition with $\sigma$ defines a (co-)limit preserving endofunctor of the identity on $\PreStk_k$, denoted by the same letter. The category of \emph{perfect prestacks}\footnote{Recall that we fixed a regular cardinal $\kappa$ in \refsect{DM.prestacks.recollections}.} is the \ii-category
\[
\PreStk_k^\perf\defined\Fun\big((\AffSch_k^\perf)^\opp,\Ani\big).
\]
The natural restriction functor $\on{res}\co \PreStk_k\to \PreStk_k^\perf$ clearly preserves limits and colimits.

\lemm \mylabel{perfect.basics}
There is an adjunction
\[
\incl : \;\PreStk_k^\perf \;\leftrightarrows \;\PreStk_k\; : \res
\]
which under the Yoneda embedding extends the adjunction $\AffSch_k^\perf\leftrightarrows \AffSch_k$.
Both functors $\incl$ and $\res$ preserve colimits, and satisfy $\on{res}\circ \on{incl} =\id$.
In particular, $\on{incl}$ is a full embedding, and $\sigma$ is an equivalence on objects in the essential image.
\xlemm

\pf
The adjunction $(\incl,\res)$ exists by \cite[Proposition~5.2.6.3]{Lurie:Higher}.
Here $\incl$ is the colimit preserving functor constructed in \cite[Theorem~5.1.5.6]{Lurie:Higher}.
Hence, the claim $\res \circ \incl = \id$ follows, by continuity of both functors, from the corresponding fact for affine schemes.
\xpf


\defi
\thlabel{colimit.perfection}
The
{\em (colimit) perfection} is the endofunctor on $\PreStk_k$ defined by
\[
(\str)^\perf\defined \on{incl}\circ \on{res}(\str).
\]
For $X\in\PreStk_k$, $X=\colim_{T\to X} T$ with $T\in\AffSch_k$, it is computed as $X^\perf=\colim_{T\to X} \lim_\sigma T$.
\xdefi

In the following, we identify $\PreStk_k^\perf\simeq \on{incl}(\PreStk_k^\perf)\subset \PreStk_k$. This is the full subcategory of all objects $X$ such that the counit of the adjunction $X^\perf\to X$ is an equivalence.
Note that the perfection of ind-objects produces ind-(perfect objects), so that the functor $(\str)^\perf$ is well-suited for working with ind-schemes.

The following lemma appeared in an early arxiv version of \cite{RicharzScholbach:Intersection}.
A more refined statement for $\SH[\frac 1 p]$, the stable $\A^1$-homotopy category localized at the characteristic (as opposed to $\DM$ with rational coefficients) is due to Elmanto and Khan \cite{ElmantoKhan:Perfection}.

\lemm
\mylabel{Frobenius.macht.nix} 
Let $X\in\PreStk_k$ be any prestack.
\begin{enumerate}
\item The pullback $\sigma^!\co \DM(X)\to \DM(X)$ is equivalent to the identity.
\item The counit of the adjunction $\pi\co X^\perf\to X$ induces an equivalence $\pi^!\co \DM(X)\simeq \DM(X^\perf)$.
\end{enumerate}
\xlemm
\pf
For i), first let $X\in \AffSch_k^\ft$. Invariance of $\DM$ under pullback along finite surjective radicial maps \cite[Proposition~2.1.9]{CisinskiDeglise:Triangulated} shows that the natural map $\id \r \sigma_*$ is invertible.
Thus, we get $\id \simeq \sigma^!$ by adjunction (using $\sigma_*=\sigma_!$). This equivalence extends to $\PreStk_k$ by continuity.

For ii), we reduce to the case $X\in\AffSch_k$ by continuity.
Then $X^\perf=\lim_\sigma X=\lim_{X\to T}\lim_\sigma T$ for $T\in \AffSch_k^\ft$, and we reduce further to the case $X=T\in \AffSch_k^\ft$.
In this case, $\pi^!\co \DM(X)\to \DM(X^\perf)=\colim_{\sigma^!}\DM(X)$ is the natural map which is an equivalence by i).
\xpf

\defi
\mylabel{limit.perfection}
The {\em \textup{(}limit\textup{)} perfection} is the endofunctor on $\PreStk_k$ defined by
\[
\lim_\sigma X \defined \lim \big(\dots \overset{\sigma}{\r} X\overset{\sigma}{\r} X\overset{\sigma}{\r} X\big),
\]
where the limit is indexed by the positive integers.
\xdefi

Restricted to schemes, this functor is used in \cite{Zhu:Affine, BhattScholze:Projectivity}. 
It gives another way of constructing prestacks on which $\sigma$ acts as an equivalence.
In general, the natural map 
$$X^\perf\to \lim_\sigma X$$ 
is not an equivalence.
For a counter-example, take the ind-scheme $\A^\infty = \colim (\A^0 \stackrel 0 \r \A^1 \stackrel{\id \x 0} \r \A^2 \dots)$. For a $k$-algebra $R$, the natural map
\[
(\bbA^\infty)^\perf(R)=\oplus_{i\geq 0} \big(\lim_\sigma R\big)\to \lim_\sigma \big(\oplus_{i\geq 0}R\big)=\lim_\sigma \bbA^\infty(R).
\]
is bijective if $R$ is reduced (so that Frobenius is injective), but not in general (e.g., take $R=k[\varpi^{1\over{p^\infty}}]/(\varpi)$). However, as $\DM$ is invariant under Nil-thickenings, the categories of motives on these two prestacks are equivalent:

\coro 
\mylabel{tomaten.auf.den.augen}
For any prestack $X\in\PreStk_k$, there is a commutative diagram of prestacks
\[
\xymatrix{(\lim_\sigma X)^\perf \ar[d] \ar[r]_{\phantom{hh}\simeq}^{\phantom{hhhh}\pi^\perf} & X^\perf \ar[d]\\
\lim_\sigma X \ar[r]^{\pi} & X,
}
\]
where $\pi$ is the natural projection. In particular by \myref{Frobenius.macht.nix} ii\textup{)}, the $!$-pullback induces equivalences $\DM(X)\simeq \DM(\lim_\sigma X)\simeq \DM(X^\perf)$.
\xcoro
\pf
We have to show that $\pi^\perf$ is an equivalence, and we claim that $\res\circ \pi$ is already an equivalence. Namely, $\Hom(T,\lim_\sigma X)\simeq \Hom(T,X)$ for any $T\in \AffSch_k^\perf$. Indeed, for any map $f\co T\r X$ and any $r\geq 0$, the diagram
\[
\xymatrix{T \ar[d]_{\sigma^r}^{\simeq} \ar[r]^{f} & X\ar[d]^{\sigma^r}\\
T \ar[r]^{f} & X,
}
\]
commutes up to equivalence.
\xpf

\rema
\mylabel{Zariski_sheaf}
If $X$ is a scheme, then the canonical map $X^\perf\to \lim_\sigma X$ induces an isomorphism after Zariski sheafification.
Indeed, if $U\to X$ is an open affine subscheme, then $U^\perf=\lim_\sigma U\to \lim_\sigma X$ is an open affine subscheme by \cite[Lemma~3.4 (ix)]{BhattScholze:Projectivity}. 
This easily implies $(X^\perf)^\Zar\cong\lim_\sigma X$. 
\nts{To check this equivalence, it is enough to compare $R$-points of $X^\perf$ with those of $\lim_\sigma X$, where $R$ is a local ring.
The closed point in $\Spec R$ maps to some $x \in U^\perf$, where $U$ runs over the open subschemes of $X$ defining a covering. The induced map $\Spec R \r \Spec \calO_{U^\perf, x$}$ shows that the image of $\Spec R$ lies in $U^\perf = \lim_\sigma U$.
Hence our map is surjective.
For the injectivity suppose $x, y: \Spec \R \r X^\perf$ are two points with $R$ local such that they get mapped to the same point in $\lim_\sigma X$. 
Since nothing happens on the topological level, we may assume $x$ and $y$ have the same set-theoretic image; by the preceding argument they hence factor over the same $U^\perf$ where $U$ is some affine part of a covering of $X$.
Then, use that $\lim_\sigma U \r \lim_\sigma X$ is a monomorphism (being a limit of such), hence the two points agree in $\lim_\sigma U = U^\perf.$}
\xrema

\subsection{Motives on perfect schemes}
\label{sect--motives.perfect.schemes}

In order to obtain a six functor formalism, we need to put suitable finiteness assumptions on the objects. 

\defi
\mylabel{Sch.pfp}
(\cite[Definition~A.13]{Zhu:Affine}, \cite[Proposition~3.11]{BhattScholze:Projectivity})
A map $f\co X \r Y$ of perfect qcqs schemes is called \emph{perfectly of finite presentation} or just \emph{pfp} if the induced map $\Spec A \r \Spec B$ of any open affines in $X$, resp.~$Y$ factors as $\Spec A \simeq \lim_\sigma (\Spec A_0)\r \Spec A_0 \r \Spec B$, where $A_0$ is a finitely presented $B$-algebra.
The full subcategory of $\Sch_k$ spanned by the qcqs pfp perfect $k$-schemes, or just {\em pfp schemes}, is denoted by $\Sch_k^\pfp$.
\xdefi



\lemm
\mylabel{loc.cat}
\begin{enumerate}

\item 
\label{item--pfp.ess.surj}
\textup{(}\cite[Proposition~A.17]{Zhu:Affine}\textup{)} The functor $\Sch_k^\ft\to \Sch_k^\pfp$, $X\mapsto \lim_\sigma X$ is full and essentially surjective.
Thus, each object $X\in\Sch_k^\pfp$ and similarly each morphism $f$ admits a {\em model}, i.e., a finite type $k$-scheme $X_0$ such that $X=\lim_\sigma X_0$.

\item
It induces an equivalence of categories $\Sch_k^\ft[\calS^{-1}]\cong \Sch_k^\pfp$ on the localization with respect to the class $\calS$ of finite radicial surjective maps in $\Sch_k^\ft$.

\item 
\label{item--pfp.equivalence.big.etale}
The equivalence in ii\textup{)} induces an equivalence on small \'etale sites.

\item
\label{item--pfp.localization}
If $Z \subset X \supset U$ is a diagram consisting of a closed and complementary open immersion in $\Sch_k^\pfp$, then there is a model $Z_0 \subset X_0 \supset U_0$ again consisting of a closed and open immersion.

\end{enumerate}
\xlemm

\pf
Part iii) follows from ii), see \StP{04DZ}. 
A map $f$ in $\Sch_k^\ft$ lies in $\calS$
iff it is a universal homeomorphism, cf.~\StPd{04DF}{01WJ}{01S4} or, equivalently by \cite[Lemma~3.8]{BhattScholze:Projectivity}, iff $\lim_\sigma f$ is an isomorphism. 
It remains to show that any functor $\varphi\co \Sch_k^\ft \r C$ to a category $C$ that sends maps in $\calS$ to isomorphisms factors uniquely over the perfection functor.
By Zhu's result in \refit{pfp.ess.surj} $\psi\co \Sch^\pfp_k \r C$ is unique if it exists. Also, to show the existence we may replace $\Sch^\pfp_k$ by the full subcategory of $\Sch_k^\pfp$ spanned by the image of the perfection functor. 
We define $\psi$ on objects by fixing a model $X = \lim_\sigma X_0$ for each $X$, and set $\psi (X) := \varphi(X_0)$.
By definition \cite[Proposition~3.11]{BhattScholze:Projectivity}, given any pfp map $f: X \r Y$, there is some $n \gg 0$ such that we have the diagram below, where $X_0^{(n)}$ is the $n$-th stage of the limit defining $X$ (it is abstractly isomorphic to $X_0$, but the structural map to $k$ differs by $\sigma^n$):
$$\xymatrix{
X = \lim_\sigma X_0 \ar[r]^f \ar[d] & Y = \lim_\sigma Y_0 \ar[d] \\
X_0^{(n)} \ar[d]^{\sigma^n} \ar[r]^{f_0^{(n)}} & Y_0 \\
X_0.
}$$
We define $\psi(f)$ as the composite $\varphi(f_0^{(n)}) \circ \varphi(\sigma^n)^{-1}$. This is well-defined independently of $n$.
One immediately checks that it is also compatible with composition. 
This shows the existence of a factorization, hence ii).

For \refit{pfp.localization}, let $X_0$ be a reduced model of $X$.
Then $|X|=|X_0|$ on the underlying topological spaces, and we define $Z_0$ (resp.~$U_0$) as the reduced closed (resp.~open) subscheme $X_0$ defined by the closed subset $|Z|$ (resp.~open subset $|U|$) in $|X|=|X_0|$.
The limit perfection preserves open and closed immersions by \cite[Lemma~3.4]{BhattScholze:Projectivity}, so that $Z_0 \subset X_0 \supset U_0$ is a model for $Z \subset X \supset U$.
\xpf

\theo
\mylabel{DM.perfect.schemes}
Motives on pfp schemes enjoy a six-functor formalism with the following properties.
Throughout, let $X$, $Y$ denote objects, resp.~$f$ a morphism in $\Sch^\pfp_k$. 
Let $X_0$, $Y_0$, resp.~$f_0$ be models thereof.

\begin{enumerate}
\item
\label{item--DM.pfp}
The presheaves $\DM^*$ and $\DM^!$ on $\Sch^\ft_k$ \textup{(}given on objects by $X \mapsto \DM(X)$, and on morphisms by $f \mapsto f^*$, resp.~$f^!$\textup{)} factor uniquely over functors
$$\DM^*, \DM^!\co (\Sch^\pfp_k)^\opp \to \DGCat_\cont.$$
On objects $X \in \Sch^\pfp_k$, these two presheaves take the same values, denoted simply $\DM(X)$.
The functor $\DM^!$ agrees with the restriction of the presheaf $\DM$ in \myref{DM.prestacks.construction}.\refit{DM!.construction} from $\PreStk_k$ to $\Sch^\pfp_k$.

\item
\label{item--DM.pfp.prestacks}
The natural map  $p \co X \r X_0$ induces an equivalence of \ii-categories
$$p^! \co \DM(X_0) \stackrel \cong \lr \DM(X), \eqlabel{DM.perfection}$$
which preserves the monoidal unit 1 and is compatible with all functors below.
In particular, if $X$ has a regular model $X_0 / k$, then 
$$\Hom_{\DM(X)}(1_X, 1_X(n)[m])\,=\,(K_{2n-m}(X_0) \t \Q)^{(n)} = \CH^n(X_0, 2n-m)_\Q,$$
where the terms at the right are the $n$-th Adams eigenspace in the algebraic $K$-theory of $X$ and Bloch's higher Chow group, both tensored with $\Q$.

\item
\label{item--functoriality.pfp}
There are adjunctions
$$f^* : \DM(Y) \rightleftarrows \DM(X) : f_*.$$
$$f_! : \DM(X) \rightleftarrows \DM(Y): f^!.$$
If $f$ has a smooth model, there is an adjunction
$$f_\sharp : \DM(X) \rightleftarrows \DM(Y): f^*.$$
These functors are functorial in $f$.
Under the equivalence \textup{\refeq{DM.perfection}}, they correspond to $(f_0)^*$ etc.

\item
\label{item--purity}
If $f$ has a proper model, then $f_! = f_*$; if $f$ has an étale model, then $f^! = f^*$. More generally, if $f$ has a smooth model $f_0$ of relative dimension $n$, there is a functorial equivalence \textup{(}called relative purity\textup{)}
$$f^! \simeq f^* (n)[2 n].\eqlabel{relative.purity.perfect}$$
Yet more generally, the latter equivalence holds if $f$ is a perfectly smooth map of relative dimension $n$ in the sense of \cite[Definition~A.25]{Zhu:Affine}.

\item
\label{item--Gm.pfp}
For the projection $p\co (\Gm)^\perf \x_k X \r X$, and any $M \in \DM(X)$, the map $p_\sharp p^* M[-1] \r M[-1]$ in $\DM(X)$ is a split monomorphism.
The complementary summand is denoted by $M(1)$.
The functor $M \mapsto M(1)$ is an equivalence with inverse denoted by $M \mapsto M(-1)$.

\item
\label{item--syno.compact}
The category $\DM(X)$ is compactly generated by the objects $t_\sharp 1(n)$, $t\co T \r X$ pfp smooth and $n \in \Z$. In particular, the monoidal unit $1_X \in \DM(X)$ is compact.
The functors $f_\sharp$, $f_*$, $f^*$, $f_!$, and $f^!$ preserve compact objects and preserve arbitrary \textup{(}homotopy\textup{)} colimits.

\item
\label{item--syno.projection.formula}
For any pfp map $f$, there is a \emph{projection formula} $(f_! M) \t N \simeq f_! (M \t f^* N)$.

\item
\label{item--syno.duality}
If $p\co X\r \Spec k$ denotes the structural map, the \emph{dualising functor}
$$\Du_X \defined \IHom(-, p^! 1)\eqlabel{define.duality.functor}$$
is a contravariant involution on the subcategory $\DM(X)^\comp$ of compact objects, i.e., $\Du_X \circ \Du_X = \id$. Furthermore, on compact objects, there are equivalences
$$\Du_Y f_! = f_* \Du_X,\;\; f^* \Du_Y = \Du_X f^!.$$

\item
\label{item--syno.localization}
For a pfp closed immersion $i\co Z \r X$ with open complement $j\co U \r X$, 
there are equivalences
$$\DM(X) = \laxlim \left (\DM(U) \stackrel{i^! j_!} \lr \DM(Z) \right),$$
$$\DM(X) = \laxlim \left (\DM(U) \stackrel{i^* j_*} \lr \DM(Z) \right).$$
Here $\laxlim$ denotes the lax limit \cite{GepnerHaugsengNikolaus:Lax} whose objects can be colloquially described as triples $(M \in \DM(U), N \in \DM(Z), \alpha: i^! j_! M \r N)$, where $\alpha$ is any map \textup{(}and the non-lax limit would be the subcategory where this map is an isomorphism\textup{)}. 
Equivalently, the \textup{(}co\textup{)}units of the adjunctions above assemble into so-called \emph{localization} homotopy fiber sequences
$$i_! i^! \r \id \r j_* j^* \stackrel {[1]\;} \r,\eqlabel{localizationOne.pfp}$$
$$j_! j^! \r \id \r i_* i^*\stackrel {[1]\;} \r.\eqlabel{localization.pfp}$$

\item
\label{item--syno.base.change}
For a cartesian diagram in $\Sch_S^\pfp$ \textup{(}necessarily consisting of pfp maps\textup{)}
$$\xymatrix{
X' \ar[r]^{g'} \ar[d]^{f'} & X \ar[d]^f \\
Y' \ar[r]^g & Y,
}$$
there are natural equivalences
\begin{align}
g^! f_* & \stackrel \cong \lr f'_* g'^!,\label{eqn--base.change.pfp}\\
f^* g_! & \stackrel \cong \lr g'_! f'^*.
\end{align}

\item
\label{item--syno.homotopy}
The category $\DM$ is \emph{perfectly homotopy-invariant} in the sense that for the projection map $p\co (\A^n)^\perf \x_k X \r X$ for any $n\in\Z_{\geq 0}$, the counit and unit maps $p_\sharp p^* \r \id$ and $\id \r p_* p^*$ are functorial equivalences in $\DM(X)$ \cite[2.1.3]{CisinskiDeglise:Triangulated}.

\item
\label{item--syno.descent}
The presheaves $\DM^*$ and $\DM^!$ are sheaves for the h-topology on pfp schemes, i.e., if $f\co X \r Y$ admits a model which is an $h$-covering in $\Sch_k^\ft$, then
the natural map
$$\DM(Y) \;\r\; \lim \left ( \DM(X) \stackrel[(p_2)^*]{(p_1)^*}\rightrightarrows \DM(X \x_Y X) \dots \right )$$
is an equivalence, and likewise with !-pullbacks instead.
We refer to this property as \emph{h-descent}.

\item
\label{item--weight.pfp}
Suppose $X$ is separated.
The category $\DM(X)$ is equipped with a \emph{weight structure} 
\[
(\DM(X)^{\wt \le 0}, \DM(X)^{\wt \ge 0}).
\]
If $X$ admits a regular model, $1_X$ is in the \emph{heart} $\DM(X)^{\wt = 0} = \DM(X)^{\wt \le 0} \cap \DM(X)^{\wt \ge 0}$ of this weight structure.
Moreover, $f^*$ and $f_!$ are weight-left exact \textup{(}preserve ``$\wt \le 0$''\textup{)} while $f_*$ and $f^!$ are weight-right exact \textup{(}preserve ``$\wt \ge 0$''\textup{)}. 

\item
\label{item--realization}
Let $\ell$ be a prime number different from $p = \Char k$.
There is a \emph{$\ell$-adic realization functor}
$$\rho_\ell\co \DM(X) \r \D_\et(X, \Ql) := \Ind(\D^\bound_\cstr(X, \Ql))$$
taking values in the ind-completion of the bounded derived category of constructible $\Ql$-adic étale sheaves constructed in \cite[§A.3]{Zhu:Affine}.
It commutes with the six functors $f_*, f^*, f_!, f^!, \t, \IHom$.
\end{enumerate}
\xtheo

\pf
\refit{DM.pfp} follows from \myref{loc.cat} and the above-mentioned fact that $f^*$ and $f^!$ are equivalences for any finite radicial surjective map $f$ \cite[Proposition~2.1.9]{CisinskiDeglise:Triangulated}.

\refit{DM.pfp.prestacks}: 
The prestack approach to $\DM(X)$ yields $\DM(X) \simeq \DM(X_0)$ by \myref{tomaten.auf.den.augen}.

With \refit{DM.pfp} and \refit{DM.pfp.prestacks} in hand, the remaining statements mostly follow directly from the classical assertions for finite-type $k$-schemes, which for the large majority are due to Ayoub and Cisinski--Déglise and, for functoriality for non-separated maps and h-descent proven in \cite[Theorem~\iref{theo--DM.descent}, Proposition~\iref{prop--unseparated}]{RicharzScholbach:Intersection}.
We will use these classical properties below without further comment, and refer to \irefsy{motives} for a detailed list of references.

The equivalence \refeq{DM.perfection} is immediate from \myref{loc.cat} using invariance of $\DM$ under finite radicial surjective maps (hence $\DM$ is well-defined on the localization) and \myref{tomaten.auf.den.augen}.
The statements in \refit{functoriality.pfp}--\refit{syno.compact}, \refit{syno.base.change}, \refit{syno.homotopy} are then direct consequences of the corresponding assertions for finite type $k$-schemes.

For \refit{syno.projection.formula} note that the functoriality $f_!$ and $f^*$ and also the tensor product are induced, by definition, from $(f_0)_!$, $(f_0)_*$ and the tensor product on $\DM(X_0)$ etc.
Similarly for the internal Hom in \refit{syno.duality}.

The purity equivalence in the case of a perfectly smooth map $f$ follows from the smooth case and the étale descent property in \refit{syno.descent}.

In the definition of $\Du_X$ as stated in \refeq{define.duality.functor}, note that the pullback functor $p^!$ exists by virtue of the general formalism of motives on prestacks. 
Under the equivalence $\DM(X) = \DM(X_0)$, it can also be computed as $\IHom(-, p_0^! 1)$, where $p_0 : X_0 \r k$ is the structural map. This again reduces all the assertions to the classical duality statements for motives on finite type $k$-schemes.

In the situation of \refit{syno.localization}, apply \myref{loc.cat}.\refit{pfp.localization} to get the localization fiber sequences \refeq{localizationOne.pfp}, \refeq{localization.pfp} from the corresponding statements for $\DM$ on $\Sch_k^\ft$ due to Ayoub.
In the presence of the base change equivalence in \refit{syno.base.change},
these fiber sequences are equivalent to the formulation involving the lax limit, cf.~the argument in \cite[Example~4.1.6]{ArinkinGaitsgory:Category}.

For \refit{weight.pfp} note that $X$ is separated iff it admits a separated model $X_0$, so that the weight structure on $\DM(X_0)$ can be used.

In the same vein, \refit{realization} follows from the equivalence \refeq{DM.perfection} and its $\ell$-adic counterpart $\D^{\bound}_\cstr(X) = \D^{\bound}_\cstr(X_0)$ \cite[§A.3.1]{Zhu:Affine}.
\xpf

\subsection{Motives on perfect ind-schemes} 
As an extension of \myref{Sch.pfp}, we consider the full subcategory $\IndSch_k^\pfp$ of ordinary presheaves consisting of strict $\aleph_0$-ind-(pfp schemes), i.e., objects are {\em countable} filtered colimits of pfp schemes with transition maps being closed immersions.
See for example \irefsect{ind.schemes}, or the notes \cite[\S1]{Richarz:Notes} for details on ind-schemes.

For brevity, the objects of $\IndSch_k^\pfp$ are called {\em ind-(pfp schemes)}. 
Similarly to \cite[Lemma~1.10]{Richarz:Notes} one shows that this category has final object $\Spec(k)$, is closed under fibre product and countable disjoint unions.
As an example we note that if $X=\colim X_i$ is an ind-scheme of ind-(finite type) over $k$, then the Zariski sheafification of the perfection $(X^\perf)^\Zar=\colim \lim_\sigma X_i$ is an ind-(pfp scheme) by \myref{Zariski_sheaf}.
We have the following basic lemma.

\lemm
\mylabel{basic.ind.pfp}
Every quasi-compact map between ind-\textup{(}pfp schemes\textup{)} is schematic and pfp.
\xlemm
\pf
Let $X\to Y$ be a quasi-compact map of ind-(pfp schemes). 
Choosing a presentation $Y=\on{colim} Y_i$ by pfp schemes, it is enough to show that each $X\x_YY_i$ is a pfp scheme.
We reduce to the case where $Y=Y_i$ is a pfp scheme, and in particular quasi-compact.
By assumption $X=\colim X_i$ is quasi-compact as well.
Being an object in $\IndSch_k^\pfp$ the underlying topological space $|X|=\colim |X_i|$ is Jacobson.
Hence, the sequence of reduced schemes $X_{i}=X_{i,\red}$ stabilizes by \cite[Corollary~1.24]{Richarz:Notes} which shows that $X=X_i$, $i>\!\!>0$ is a pfp scheme.
\xpf

\theo 
\mylabel{functor.adic}
Motives on ind-\textup{(}pfp schemes\textup{)} $X \in \IndSch_k^\pfp$ satisfy the properties 
listed in \myref{DM.perfect.schemes}, with the following adjustments:

\begin{itemize}
\item
The statement in \refit{DM.pfp.prestacks} has no analogue for ind-\textup{(}pfp schemes\textup{)}, since they need not have models among ind-schemes of ind-\textup{(}finite type\textup{)}.

\item
The description of compact generators \textup{(}cf.~\refit{syno.compact} above\textup{)} is as follows: the category $\DM(X)$ is compactly generated by the images of the compact objects in $\DM(X_i)$, where $X = \colim X_i$ is a presentation.

\item
The weight structure on $\DM(X)$ exists if there is a presentation $X = \colim X_i$, where the $X_i \in \Sch_k^\pfp$ are separated. 

\item
The functor $f^*$ exists \textup{(}and satisfies the properties as stated in the remaining items\textup{)} if $f$ is quasi-compact \textup{(}hence schematic and pfp\textup{)}. 

\item
If $X$ is componentwise quasi-compact, then $p^* 1$ is a monoidal unit, for the structural map $p\co X \r \Spec k$.
In general, $\DM(X)$ only carries a non-unital symmetric monoidal structure.

\item
The condition that $f$ has a smooth model \textup{(}resp.~is a pfp closed immersion, pfp open immersion\textup{)} has to be replaced with ``$f$ has a schematic smooth model $f_0$'' \textup{(}has a model which is a schematic pfp closed, resp. open immersion\textup{)}.
The h-covering in \refit{syno.descent} has to be schematic.
The term ``proper'' has to be replaced by ``ind-\textup{(}perfectly proper\textup{)}'' in \refit{purity}.
\end{itemize}
\xtheo

\pf 
Recall from \irefth{motives.Ind-schemes} that the corresponding assertions hold for motives on ind-schemes.
From there, the theorem as stated follows from \myref{loc.cat}, in the same way as the assertions for motives on pfp schemes follow from those on finite-type $k$-schemes.
\xpf

\section{Stratified Tate motives for perfect ind-schemes}

With this convenient extension of the six functor formalism for the categories $\DM$ of motives on pfp schemes and their ind-variants, we can quickly adapt the contents of \cite[§§\iref{sect--DTM}--\iref{sect--DTM.Fl}]{RicharzScholbach:Intersection}.

\defi
A \emph{stratified ind-\textup{(}pfp scheme\textup{)}} is a map of ind-(pfp schemes)
\[
\iota\co X^+ := \bigsqcup_{w \in W} X_w \r X,
\] 
where $\iota$ is bijective on the level of underlying sets, each $X_w$ is a pfp scheme, each restriction $\iota|_{X_w}$ is an immersion and the topological closure of each $\iota(X_w)$ is a union of strata. 
Here $W$ is a countable index set which is partially ordered by the closure relations on the strata, i.e., $\iota(X_v)$ lies in the closure of $\iota(X_w)$ if and only if $v\leq w$.

Since each $X_w$ is pfp (hence quasi-compact), each restriction $\iota|_{X_w}\co X_w\to X$ is a quasi-compact map of ind-(pfp schemes) and hence schematic and pfp by \myref{basic.ind.pfp}.
A map of stratified ind-(pfp schemes) is required to be quasi-compact (hence, again, schematic and pfp), and to map strata to strata.
\xdefi

\defi 
\mylabel{strat.adic}
A {\em perfect cell} is a $k$-scheme isomorphic to the perfection of a cell $\A^n_k \x \GmX k^r$.  
A pfp scheme $X$ is {\em perfectly cellular} if it admits a smooth model over $k$, and admits a stratification into perfect cells. 
(We do not require that the smooth model admits a stratification into cells.)
A {\em perfectly cellular stratified}, or just {\em pcs ind-scheme} is a stratified ind-(pfp scheme) $X$ where each $X_w$, $w\in W$ is a perfectly cellular $k$-scheme.
\xdefi

For a pfp scheme $X$, or a disjoint union of those, the category of Tate motives $\DTM(X)\subset \DM(X)$ is the full subcategory generated under arbitrary shifts and colimits by $1_X(n)$, $n \in \Z$.
By \myref{DM.perfect.schemes} \refit{DM.pfp.prestacks}, it is equivalent to $\DTM(X_0)$ as defined in \irefde{Tate.geometry}, where $X_0$ is any model of $X$. 
Indeed, the equivalence $\sigma^!\co \DM(X_0) \r \DM(X_0)$ restricts to an equivalence of the $\DTM$-subcategories by \myref{DM.perfect.schemes}~\refit{DM.pfp.prestacks}.

For the next definition, we note that by quasi-compactness of $\iota$ the pullback $\iota^*$ is well-defined, cf.~\myref{functor.adic}.

\defi
\mylabel{Whitney-Tate.pfp}
A stratified ind-(pfp scheme) $X$ is \emph{Whitney--Tate} if $\iota^* \iota_*\co \DM(X^+) \r \DM(X^+)$ preserves the subcategory $\DTM(X^+)$, i.e., $(\iota_w)^* (\iota_v)_*$ induces $\DTM(X_v)\to \DTM(X_w)$ for all $v,w\in W$. 
In this case, the category of \emph{stratified Tate motives}
$$\DTM(X, X^+) \defined \laxlim_{w \in W} \DTM(X_w)$$
is well-defined. 
Here the transition functors used to form the lax limit is $(\iota_w)^* (\iota_v)_*$.
Since $\DM(X)$ is equivalent to the lax limit of the categories $\DM(X_w)$, the category $\DTM(X, X^+)$ is a full subcategory of $\DM(X)$.
A map $\pi$ between stratified ind-(pfp schemes) is a \emph{Whitney--Tate map} if $\pi_*$ preserves stratified Tate motives.
\xdefi

The characterization of Whitney--Tate stratifications in \irefdele{Whitney.Tate.condition} holds literally. 
Also the perfection of a cellular stratified $k$-scheme of finite type in the sense of \irefde{cellular} is pcs.
Conversely, every pcs scheme $\sqcup X_w\to X$ admits a model $\sqcup X_{w,0}\to X_0$ where each $X_{w,0}$ is smooth.
We now show that the $t$-structures on the individual strata $X_w$, which are induced from a model $X_{w,0}$, glue together for any pcs (ind-)scheme.

Recall from \cite{Levine:Tate} that the \emph{Beilinson--Soul\'e vanishing conjecture} for $k$ is equivalent to the existence of a $t$-structure on $\DTM(\Spec k)$ such that $1_{\Spec k}(r)$ is in the heart for all $r \in \Z$.
By Quillen's work, this conjecture holds for $k = \Fq$.
By Harder's work, it holds for a function field $\Fq(t)$ and therefore, by the compatibility of algebraic $K$-theory with filtered colimits of rings, it also holds for its perfection $k=\Fq(t)^\perf$.

\defilemm
\mylabel{DTM.t-structure.pfp}
Suppose $k$ satisfies the Beilinson--Soul\'e vanishing conjecture.
Let $X$ be a pcs Whitney--Tate ind-scheme over $k$.
Then each $\DTM(X_w)$ carries a unique $t$-structure such that $1_{X_w}(n)[d_w]$, $d_w:=\dim(X_w)$ is in the heart.

The subcategories 
$$\eqalign{
\DTM(X, X^+)^{\le 0} & := \big\{M \in \DTM(X) \;|\; \iota^* M \in \DTM(X^+)^{\le 0} \big\} \cr
\DTM(X, X^+)^{\ge 0} & := \big\{M \in \DTM(X) \;|\; \iota^! M \in \DTM(X^+)^{\le 0} \big\}
}$$
define a compactly generated $t$-structure on $\DTM(X, X^+)$.
The heart of this $t$-structure is denoted by $\MTM(X, X^+)$ or just by $\MTM(X)$ if the stratification is clear from the context.
\xdefilemm

\pf
By assumption there exists a finite stratification $ \sqcup_iX_w^{(i)}\to X_w$ into perfect cells.
We may choose a model $\sqcup_i X_{w,0}^{(i)}\to X_{w,0}$ where both source and target are smooth.
(Here the schemes $X_{w,0}^{(i)}$ are cells, but this map is only a stratification after perfection.)
Then the first claim follows from \irefle{t-structure.stratum} using that 
the equivalence $\DTM(X_w) = \DTM(X_{w,0})$ in \myref{DM.perfect.schemes} \ref{item--DM.pfp.prestacks} commutes with the six functors (by construction).
The second statement is also shown similarly as in \emph{op.~cit.}: for any $w \in W$, the $t$-structure on $\DTM(X_{\leq w}, X^+)^\comp$ arises by glueing the $t$-structures on the strata $X_v$, $v \le w$.
Given that the compact objects in $\DM(X,X^+)$ are supported on some finite union of $X_{\leq w}$, this gives a $t$-structure on $\DTM(X,X^+)^\comp$ and therefore a compactly generated $t$-structure on $\DTM(X,X^+)$. 
\xpf

Restricting the $\ell$-adic realization functor to Tate motives continues to have the pleasant properties ---such as conservativity and faithfulness--- we know from the case of cellular Whitney--Tate stratified ind-schemes.
The proof of the following statement is a verbatim copy of \irefle{Tate.conservative} and \sref{realization.functor.coro}.

\prop
\mylabel{realization.pfp}
Let $X$ be a pcs Whitney--Tate ind-scheme over $k$ and suppose that $k$ satisfies the Beilinson--Soulé vanishing conjecture.
The restriction of the $\ell$-adic realization functor $\rho_\ell$ from $\DM(X)$ to $\DTM(X) := \DTM(X, X^+)$ has the following properties: 

\begin{enumerate}
\item 
it is conservative,

\item 
it creates the motivic t-structure in the sense that $M \in \DTM(X)^\comp$ is in positive (resp.~negative) degrees of the motivic t-structure iff $\rho_\ell(M)$ has the corresponding property for the perverse t-structure.
In particular, it restricts to a functor
$$\rho_\ell\co \MTM(X)^\comp \r \Perv(X, \Ql).\eqlabel{rho.ell.MTM}$$
Furthermore, this latter functor is faithful.
\end{enumerate}
\xprop

\section{Motives on Witt vector affine flag varieties}
\label{sect--Witt.vector.flag.varieties}
Let $K/\bbQ_p$ be a finite extension with uniformizer $\varpi\in \calO_K$ and residue field $k/\bbF_p$.
Let $G$ be a split reductive group over $\calO_K$. 
In this section, we introduce Tate motives on the {\em Witt vector affine flag varieties} for $G$ and prove basic properties thereof. 

\subsection{Loop groups and their affine flag varieties} 
For any $k$-algebra $R$, we denote by $W_{\calO_K}(R)$ the ring of ramified Witt vectors, cf.~e.g.~\cite[\S1.2]{FarguesFontaine:Courbe}. 
If $R$ is perfect, then
\[
W_{\calO_K}(R)=W(R)\otimes_{W(k)}\calO_K=\Big\{\sum_{i\geq 0}[\la_i]\varpi^i\;\big|\;\la_i\in R\Big\},
\]
where $W(k)\to \calO_K$ is the unique map, and where $[\cdot]\co R\to W(R)$ is the Teichm\"uller lift. 
The {\em (Witt vector) loop group} is the group-valued functor on the category of perfect $k$-algebras given by
\[
LG(R)\defined G\Big(W_{\calO_K}(R)\big[{\textstyle{1\over p}}\big]\Big).
\]
For any smooth, affine $\calO_K$-group scheme $\calG$, the \emph{Witt vector positive loop group} is the group-valued functor on the category of perfect $k$-algebras given by
\[
L^+\calG(R)\defined \calG\big(W_{\calO_K}(R)\big).
\]
If $\calG\otimes K=G\otimes K$, then $L^+\calG\subset LG$ defines a closed subgroup functor, cf.~\cite[Lemma~1.2 (i)]{Zhu:Affine}. 
Being presheaves on $\AffSch^\perf_k$, $LG$ and $L^+ \calG$ are also perfect prestacks, so that we can consider categories such as $\DM(LG)$ etc.

Fix a Borel pair $T\subset B$ in $G$ over $\calO_K$. 
Denote by $\scrA=\scrA(G,T,K)$ the associated apartment of the Bruhat--Tits building of $G(K)$ equipped with the base point $0$ corresponding to $G/\calO_K$.
The choice of $B$ induces a unique alcove $\bba_0\subset \scrA$ in the dominant chamber whose closure contains $0$.
For each facet $\bbf\subset \scrA$, there is a parahoric group scheme $\calG_\bbf$ over $\calO_K$ with generic fibre $G\otimes K$ by \cite[\S5.2]{BruhatTits:Groups2}.
Recall that $\calG_\bbf$ is a smooth, affine $\calO_K$-group scheme with connected fibers such that $\calG_\bbf(\calO_K)\subset G(K)$ fixes $\bbf$ pointwise.

\lemm 
\mylabel{p.adic.rep}
\begin{enumerate}
\item The functor $LG$ is representable by a group ind-\textup{(}perfect affine scheme\textup{)}.
\item For every smooth, affine $\calO_K$-group scheme $\calG$, the functor $L^+\calG=\lim_{i\geq 0}\calG_{i}$ is representable by the perfection of an inverse system of smooth, affine group schemes with smooth, surjective transition maps. 
For each $i\geq 0$, the kernel $\ker(\calG_{i+1}\to \calG_{i})$ is the perfection of a vector group. 
In particular, every $L^+\calG$-torsor in the fpqc-topology admits sections \'etale locally. 
If $\calG=\calG_\bbf$ for some facet $\bbf\subset\scrA$, then each $\calG_{i}=\calG_{\bbf,i}$ is moreover the perfection of a cellular $k$-scheme, i.e., the perfection of a smooth \textup{(}finite type\textup{)} $k$-scheme which can be stratified by pieces isomorphic to $\A^r_k \x \GmX k^s$.
\end{enumerate}
\xlemm
\pf For the representability statements see \cite[Proposition~9.2]{BhattScholze:Projectivity}, and the references cited therein.
For the \'etale local triviality, we note that $L^+\calG$ is already representable as a functor on the category of all $k$-algebras by the same formula by Greenberg \cite{Greenberg:Witt} (cf.~also \cite[\S1.1.1]{Zhu:Affine}, \cite[Remark~9.3]{BhattScholze:Projectivity}). 
It is the perfection of the pro-algebraic group $\lim_{i\geq 0}\calH_i$ where $\calH_i\co R\mapsto \calG(W_{\calO_K}(R)/V_\varpi^iW_{\calO_K}(R))$ and $V_\varpi$ denotes the Verschiebung on the ramified Witt vectors. 
Consequently, $L^+\calG=\lim_{i\geq 0}\calG_i$ where $\calG_i=\calH_i^\perf$. For any $k$-algebra $R$, there is a short exact sequence
\[
0\to R\to W_{\calO_K}(R)/V_\varpi^{i+1}W_{\calO_K}(R)\to W_{\calO_K}(R)/V_\varpi^iW_{\calO_K}(R)\to 0.
\]
By \irefpr{vector.extension} we have $\bbV(\calE_i)=\ker(\calH_{i+1}\to \calH_{i})$ where $\calE_i$ is the vector space given by the formula \inter{(\iref{vector_bundle})} applied with $X=\Spec(W_{\calO_K}(R))$, $\hat{D}=\colim_{i\geq 0}D_i$ for $D_i=\Spec(W_{\calO_K}(R)/V_\varpi^i W_{\calO_K}(R))$. 
Restricting back to perfect $k$-algebras, we see that $\bbV(\calE_i)^\perf=\ker(\calG_{i+1}\to \calG_{i})$ is the perfection of a vector group. Now by \cite[Lemma~A.9]{Zhu:Affine} a torsor under $\calG_i=\calH_i^\perf$ is trivial if and only if its pushout along $\calG_i\to \calH_i$ is trivial. 
In particular, each $\calG_i$-torsor is \'etale locally trivial, and each affine $\ker(\calG_{i+1}\to \calG_{i})$-torsor is trivial. 
Thus, \inter{Lemma~\iref{lemm--torsor.sequence}, Corollary~\iref{coro--etale.torsor}} applies to show that each fpqc-$L^+\calG$-torsor is \'etale locally trivial. 
Finally, let $\calG=\calG_\bbf$ for some facet $\bbf$. 
By the same argument as in \irefle{affine.proj} each $\calH_{\bbf,i}$ is cellular, and thus $\calG_{\bbf,i}=\calH_{\bbf,i}^\perf$ is perfectly cellular.
\xpf

For any facet $\bbf\subset \scrA$, we denote by $\calP_\bbf:=L^+\calG_\bbf\subset LG$ the associated parahoric subgroup. The Iwahori subgroup is denoted $\calB:=\calP_{\bba_0}$, and the positive loop group is denoted $L^+G:=\calP_{\{0\}}$.

\defi
The \emph{Witt vector \textup{(}partial\textup{)} affine flag variety} is the \'etale sheaf quotient
\[
\Fl_\bbf\defined (LG/\calP_\bbf)^\et,
\]
i.e., the \'etale sheaf associated with the (ordinary) presheaf $LG/\calP_\bbf$ in $\PreStk_k^\perf$.
\xdefi

By \myref{p.adic.rep}, the \'etale sheaf $\Fl_\bbf$ is an fpqc sheaf, so the definition agrees with \cite[\S1.1.2]{Zhu:Affine} and \cite[Definition~9.4]{BhattScholze:Projectivity}. 
In particular, it is representable by an ind-scheme $\Fl_\bbf=\colim_i X_i$ where each $X_i$ is the perfection of a projective $k$-scheme~\cite[Corollary~9.6]{BhattScholze:Projectivity}.
In particular, it is a strict ind-(pfp scheme) over $k$.

Now fix two facets $\bbf,\bbf'\subset \scrA$ contained in the closure of the base alcove. The definitions and arguments in \irefsect{loop.grps} translate (almost) literally: the left action of $\calP_{\bbf'}$ on $\Fl_\bbf$ induces a presentation
\begin{equation}\label{pres.flag.adic}
\Fl_\bbf\;=\; \colim_{w}\Fl_{\bbf, \leq w},
\end{equation}
where $w$ ranges through the double quotient $W_{\bbf'}\backslash W/W_\bbf$ of the Iwahori--Weyl group by the reflection subgroups associated with the facets, defined analogously to \inter{(\iref{IW_Indentify})}. 
Here the Schubert varieties $\Fl_{\bbf, \leq w}$ are defined as the scheme-theoretic image of the orbit map $\calP_{\bbf'}\to \Fl_\bbf, p'\mapsto p'\cdot \dot{w}\cdot e$ analogously to \irefde{Schubert.scheme}. 
Then each $\Fl_{\bbf,\leq w}$ is the perfection of a projective $k$-variety by \cite[Corollary~9.6]{BhattScholze:Projectivity}. 
The \'etale sheaf image $\Fl_{\bbf,w}\subset \Fl_{\bbf, \leq w}$ of the orbit map is representable by the perfection of a smooth $k$-variety, and identifies with the \'etale sheaf quotient $\calP_{\bbf'}/\calP_{\bbf'\cup w\bbf}$ (combine the argument of \irefle{orbit.flag} with \cite[Proposition~A.32]{Zhu:Affine}.) 
Here $\calP_{\bbf'\cup w\bbf}:=L^+\calG_{\bbf'\cup w\bbf}$ where $\calG_{\bbf'\cup w\bbf}$ is the Bruhat-Tits $\calO_K$-group scheme with geometrically connected fibers determined by the subset ${\bbf'\cup w\bbf}\subset \scrA$. 

There is a stratification
\begin{equation}\label{strat.flag.adic}
\iota\co \Fl_\bbf^+:=\bigsqcup_{w\in W_{\bbf'}\backslash W/W_\bbf} 
{\Fl}{_{\bbf,w}}
\lr \Fl_\bbf,
\end{equation}
where the closure relations of the strata are given by the Bruhat order on $W_{\bbf'}\backslash W/W_\bbf$. This follows from the existence of Demazure resolutions \cite[(1.4.1)]{Zhu:Affine} using the method in \cite[Proposition~2.8]{Richarz:Schubert}. If $\bbf'=\bba_0$, we refer to this stratification as the Iwahori stratification. We summarize important properties in the following lemma.

\lemm 
\mylabel{Fl.strata.adic}
Let $\bbf,\bbf'\subset \scrA$ be contained in the closure of the base alcove.
\begin{enumerate}
\item 
The stratification \eqref{strat.flag.adic} is a pcs stratification in the sense of \myref{strat.adic}. 
\item 
\label{item--stabilizers.connected}
The stabilizers of the $\calP_{\bbf'}$-left action on $\Fl_\bbf$ are geometrically connected. 
\item 
\label{item--change.facet.properties}
If $\calP_{\bbf'}\subset \calP_\bbf$, the projection $\pi\co \Fl_{\bbf'}\to \Fl_\bbf$ is schematic, perfectly proper, and \'etale locally on the target isomorphic to the projection $\calP_{\bbf}/\calP_{\bbf'}\x \Fl_\bbf\to\Fl_\bbf$. The induced map on the {\emph Iwahori} strata $\pi^+\co \Fl_{\bbf'}^+\to \Fl_\bbf^+$ is a Whitney--Tate map and admits a section $s^+\co\Fl_{\bbf}^+\to \Fl_{\bbf'}^+$ which is an open and closed immersion.
\item In iii\textup{)}, the \'etale quotient $\calP_{\bbf}/\calP_{\bbf'}$ is the perfection of an homogenous space. 
\end{enumerate}
\xlemm
\pf Part i) follows from the discussion above, and the following observation: For each $w\in W_{\bbf'}\backslash W/W_\bbf$, we have the Iwahori stratification 
${\Fl}{_{\bbf,w}} = \bigsqcup_{v} \leftidx{^{\bba_0}}\!{\Fl}{_{\bbf,v}}$,  where $v$ runs through $W_{\bbf'}wW_\bbf/W_\bbf$, cf.~\cite[\S1.4]{Zhu:Affine}. 
Each Iwahori orbit $\leftidx{^{\bba_0}}\!{\Fl}{_{\bbf,v}}$ is isomorphic to $(\bbA^{l(v)}_k)^\perf$ where $l(v)$ is the length in $W/W_\bbf$. 
This is proven as in \inter{(\iref{Root_Groups})}, but working with perfections everywhere:
for each affine root $\al$, there is a root subgroup $U_\al\subset LG$, $U_\al\simeq (\bbA^1_k)^\perf$, and an isomorphism
\begin{equation}\label{product_map_eq}
\bigsqcap_\al U_\al\overset{\simeq}{\lr} \leftidx{^{\bba_0}}\!{\Fl}{_{\bbf,v}}, \;\;\; (u_\al)\longmapsto \big(\sqcap u_\al\big)\cdot \dot w\cdot e,
\end{equation}
where the product (taken in any order) ranges over all affine roots $\al$ such that $(w\al)|_{\bba_0}$ takes positive values and $\al|_\bbf$ takes negative values. 
Part ii) follows from the presentation ${\Fl}{_{\bbf,w}}=\calP_\bbf/\calP_{\bbf'\cup w\bbf}$ and the geometrical connectedness of $\calP_{\bbf'\cup w\bbf}$ explained above.
Part iii) is proven analogously to \irefpr{change.facet}; note that $s^+$ is open since it is an inclusion of a disjoint union of certain Iwahori strata.
For iv), let $\calP_{\bbf,0}=\calH_0^\perf$, $\calP_{\bbf',0}=(\calH_0')^{\perf}$ be the Greenberg realization as in the proof of \myref{p.adic.rep}. 
Let $\bar{\calP}_{\bbf',0}=\im(\calP_{\bbf',0}\to \calP_{\bbf,0})$, and let $\bar{\calH}_0'=\im(\calH_0'\r\calH_0)$ which is a parabolic subgroup of $\calH_0$.
Then on \'etale quotients
\[
\calP_\bbf/\calP_{\bbf'}\;\cong\; \calP_{\bbf,0}/\bar{\calP}_{\bbf',0}\;\cong\; \lim_\sigma(\calH_0/\bar{\calH}_0').
\] 
For the second isomorphism, we use that these \'etale quotients are already fpqc sheaves, and that taking perfections preserves faithfully flat maps \cite[Lemma~3.4 (xii)]{BhattScholze:Projectivity}, hence $\bar{\calP}_{\bbf',0}\cong (\bar{\calH}_0')^\perf$.
See also \cite[Theorem~6.1]{BertapelleGonzalez:Perfection}.
\xpf

\exam 
\mylabel{simple.reflection.adic}
Let $s\in \bbS$ be a simple affine reflection. 
Let $\bbf_s$ be the unique facet of maximal dimension in the closure of $\bba_0$ such that $s(\bbf_s)=\bbf_s$, i.e., $W_{\bbf_s}$ is the subgroup of $W$ generated by $s$.
We specialize \myref{Fl.strata.adic} iii) to the case $\bbf'=\bba_0$, and $\bbf=\bbf_s$, so that $\pi\co \Fl:=\Fl_{\bba_0}\to \Fl_{\bbf_s}$ is the projection from the full affine flag variety. In this case, $\pi$ has general fiber $\calP_{\bbf_s}/\calB=\lim_\sigma\P_k$ by \myref{Fl.strata.adic} iv) (cf.~also \cite[above (1.4.1)]{Zhu:Affine}). If $w\in W$, and $w=v\cdot s$ is a reduced decomposition, then
\[
(\pi^+)^{-1}\big(\Fl_{\bbf_s,v}\big)\,=\,\Fl_{v}\sqcup\Fl_{vs}.
\]
Here $\pi^+|_{\Fl_v}$ is an isomorphism, and $\pi^+|_{\Fl_{vs}}$ is isomorphic to the projection $(\bbA^1_k)^\perf\x\Fl_v\to \Fl_v$.
\xexam

\subsection{Motives on affine flag varieties}
All results from \inter{\S\S\iref{sect--WT.Fl}--\iref{sect--DTM.double}} except \irefco{equivariant.IC} (there seems to be no such $\Gm$-action on Witt vector loop groups) translate to the setting of Witt vector affine flag varieties.
For the purposes of this paper, the salient features are as follows:

\theo 
\mylabel{Fl.MTM.adic}
Let $\bbf,\bbf'\subset \scrA$ be contained in the closure of the base alcove.
\begin{enumerate}
\item 
\label{item--MTM.Fl.pfp}
The pcs stratification \eqref{strat.flag.adic} into $\calP_{\bbf'}$-orbits on $\Fl_\bbf$ is Whitney--Tate. In particular, by \myref{Whitney-Tate.pfp} and \myref{DTM.t-structure.pfp} there are well-defined categories of motives
\[
\MTM(\Fl_\bbf)\;\subset\; \DTM(\Fl_\bbf)\;\subset\; \DM(\Fl_\bbf).
\]
The $\bbQ$-linear abelian category $\MTM(\Fl_\bbf)$ is generated (by means of extensions and direct sums) by the simple objects called \emph{intersection motives}
$$\IC_w(n) := (\iota_w)_{!*} (1 (n))$$
for $n\in \bbZ$, $w\in W_{\bbf'}\backslash W/W_\bbf$. 
These objects are defined as the Tate twists of intermediate extensions along the inclusions $\iota_w : \Fl_{\bbf,w} \r \Fl_\bbf$.
(Note that the dependency of the categories $\DTM(\Fl_\bbf)$ and $\MTM(\Fl_\bbf)$ on the choice of $\bbf'$ is suppressed in the notation.)

\item
\label{item--l-adic.Fl.pfp}
Let $\ell$ be a prime number $\ne p$.
The $\ell$-adic realization functor $\rho_\ell$ (cf.~\refeq{rho.ell.MTM}) maps $\IC_w(n)$ to the $\ell$-adic intersection complexes on the Schubert varieties $\Fl_{\bbf, \leq w}$. 
A compact object $M\in \DTM(\Fl_\bbf)^\comp$ lies in $\MTM(\Fl_\bbf)^\comp$ iff $\rho_\ell(M)$ is a perverse $\ell$-adic sheaf.
\end{enumerate}
\xtheo
\pf 
The proofs in \inter{\S\S\iref{sect--WT.Fl}--\iref{sect--DTM.double}} carry over literally, as we briefly indicate.
The Whitney--Tate property is shown in three steps: 
First, for $\bbf' = \bbf = \bba_0$, the Whitney--Tate condition in \refit{MTM.Fl.pfp} is a consequence of the geometric properties listed in \myref{simple.reflection.adic} together with an induction on the length of $w$.
Second, the case $\bbf' = \bba_0$ and arbitrary $\bbf$ follows from a general criterion \irefle{Tate.proper.descent} that allows to carry over the Whitney--Tate property along the projection map $\Fl_{\bba_0} \r \Fl_\bbf$, using the geometric properties of this map listed in \myref{Fl.strata.adic}.\refit{change.facet.properties}. Third, for both $\bbf, \bbf'$ arbitrary, we use \irefpr{DTM.G/H}: a $\calP_\bbf$-equivariant motive is Tate with respect to the stratification by $\calP_\bbf$-orbits as soon as it is Tate with respect to the finer stratification by Iwahori-orbits (i.e. $\calP_{\bba_0}$-orbits).
This latter property holds by the second step.

The category $\DTM(\Fl_{\bbf,w})$ of (unstratified) Tate motives on the strata $\Fl_{\bbf,w}$ carries a $t$-structure whose heart $\MTM(\Fl_{\bbf,w})$ is an abelian category whose simple objects are of the form $1(n)[\dim \Fl_{\bbf,w}]$ (\myref{DTM.t-structure.pfp}).
These t-structures are glued together, as is classical in the context of perverse sheaves \cite[§1.4]{BeilinsonBernsteinDeligne:Faisceaux}; see also \cite[Theorem 10.3]{SoergelWendt:Perverse} for an application of this idea in the context of stratified Tate motives.

The description of the simple generators then follows from the description of the simple generators in the categories that are glued together \irefth{generators.DTM.flag}.

\refit{l-adic.Fl.pfp} is then a direct consequence of \myref{realization.pfp}.
\xpf

\theo
\mylabel{DTM.doubly.equivariant}
\begin{enumerate}
\item
\label{item--DTM.right.left.pfp}
The full subcategories of $\DM(\calP_{\bbf'}\backslash LG/\calP_\bbf)$,
$$\DM(\calP_{\bbf'}\backslash LG/\calP_\bbf) \x_{\DM(\Fl_\bbf)} \DTM(\Fl_\bbf)
\;\;\;\text{and}\;\;\;
\DM(\calP_{\bbf'}\backslash LG/\calP_\bbf) \x_{\DM(\Fl^\opp_{\bbf'})} \DTM(\Fl^\opp_{\bbf'})$$
 agree where $\Fl^\opp_{\bbf'}  := (\calP_\bbf \bsl LG)^\et$.
\textup{(}For example, the first category consists of those motives on $\calP_{\bbf'}\backslash LG/\calP_\bbf$ whose underlying, i.e., non-equivariant motive on $\Fl_\bbf$ is a stratified Tate motive with respect to the stratification by $\calP_{\bbf'}$-orbits.\textup{)}
This category is denoted $\DTM(\calP_{\bbf'}\backslash LG/\calP_\bbf)$.

\item
\label{item--MTM.forget.pfp}
For $\bbf=\bbf'$, there is a $t$-structure on $\DTM(\calP_{\bbf'}\backslash LG/\calP_\bbf)$ such that the forgetful functors to $\DTM(\Fl_\bbf)$ \textup{(}and also to $\DTM(\Fl^\opp_\bbf)$\textup{)} are $t$-exact. 
The heart of this $t$-structure is the category of mixed Tate motives on the double quotient, denoted $\MTM(\calP_{\bbf'}\backslash LG/\calP_\bbf)$.
The functor forgetting the $\calP_\bbf$-equivariance of a mixed Tate motive,
$$\MTM(\calP_\bbf\backslash LG/\calP_\bbf)\to\MTM(\Fl_\bbf),$$ 
is fully faithful, and induces a bijection of isomorphism classes of simple objects as described in \myref{Fl.MTM.adic}.\refit{MTM.Fl.pfp}.
\end{enumerate}
\xtheo

\pf
\refit{DTM.right.left.pfp} is shown by exhibiting an explicit set of generators, in the same way as in \irefth{equivariant.DTM.flag}.
\refit{MTM.forget.pfp}: The t-structure on $\DTM(\calP_{\bbf'}\backslash LG/\calP_\bbf)$ is an instance of a t-structure on a limit of a diagram of stable \ii-categories equipped with t-structures and t-exact transition functors \irefle{t.structure.limit}.
For $\bbf = \bbf'$ (but not in general, because of the shift by the dimension of $\Fl_\bbf^w$), this t-structure is identical to the one obtained by using the underlying non-equivariant motive on $\Fl_\bbf^\opp$ \irefth{equivariant.DTM.flag}.

The full faithfulness of the forgetful functor is a general consequence of \irefpr{equivariant.MTM}:
the key point is that the stabilizers of the $\calP_\bbf$-action on $\Fl_\bbf$ are connected by \myref{Fl.strata.adic}.\refit{stabilizers.connected} and that for the perfection of a connected algebraic group $H$, the pullback functor $\MTM(X) \r \MTM(X \x H)$ is fully faithful \irefle{pi*.fullyfaithful.MTM}, which implies that the forgetful functor 
\[
\MTM(X / H) = \lim \big(\MTM(X) \rightrightarrows \MTM(X \x H)\dots\big) \r \MTM(X) 
\]
is fully faithful as well.
Finally, the intersection motives $\IC_w(n)$ are $\calP_\bbf$-equivariant  by construction and therefore also generate the subcategory $\MTM(\calP_\bbf\backslash LG/\calP_\bbf)\subset \MTM(\Fl_\bbf)$ (again, by means of extensions and direct sums).
\xpf

\subsection{The convolution product}
\label{sect--convolution.product.witt}
In this section, we discuss the convolution product on the category $\MTM(\calP_\bbf \bsl LG / \calP_\bbf)$ by indicating how the corresponding arguments in the equi-characteristic situation carry over to the case considered here.

We define the convolution product by the same formula as in \sref{convolution.product}:
$$\str\star\str := m_! p^! ( \str \boxtimes \str) \colon \DM(\calP_{\bbf} \bsl LG / \calP_\bbf) \x \DM(\calP_\bbf \bsl LG / \calP_{\bbf}) \lr \DM(\calP_{\bbf} \bsl LG / \calP_{\bbf}),$$
where $p$ and $m$ are the projection and multiplication maps:
$$\xymatrix{
\calP_{\bbf} \backslash LG / \calP_\bbf \x \calP_\bbf \backslash LG / \calP_{\bbf} &
\calP_{\bbf} \backslash LG \x^{\calP_\bbf} LG / \calP_{\bbf} \ar[l]_(.45)p \ar[r]^(.55)m &
\calP_{\bbf} \backslash LG / \calP_{\bbf}}.$$
As in \irefpr{boxtimes}, the exterior product 
$$\boxtimes : \DM(\calP_\bbf \bsl \Fl_\bbf) \x \DM(\calP_\bbf \bsl \Fl_\bbf) \r \DM(\calP_\bbf \bsl \Fl_\bbf \x \calP_\bbf \bsl \Fl_\bbf)$$
is defined using that the kernels of $\calP_\bbf \r \calP_{\bbf, i}$ are perfectly split pro-unipotent (i.e., of the form $\lim G_i$, where the underlying schemes of $G_0$ and all of the $\ker (G_{i+1} \r G_i)$ are perfections of vector groups, cf.~the proof of \myref{p.adic.rep}) and the perfect homotopy invariance of $\DM$ (\myref{DM.perfect.schemes}.\refit{syno.homotopy}).
The pushforward $m_!$ exists by \myref{functor.adic} (and the descent of functoriality expressed in \irefle{functoriality.equivariant}) since the multiplication map induces a schematic map of ind-schemes $\Fl_\bbf \xtw \Fl_\bbf := (LG \x^{\calP_\bbf} LG/\calP_\bbf)^\et \r \Fl_\bbf$.
Note that the convolution of two compact objects is again compact.

\rema
The exterior product for motives on $\Fl_\bbf^\opp / \calP_\bbf$ can be used to define another convolution product functor.
At least on the level of the homotopy categories $\Ho(\DM(\calP_{\bbf} \bsl LG / \calP_\bbf))$, these two functors are isomorphic as can be shown like in \sref{convolution.HoDM.independent}, using that $LG_{\le w} = \calP_\bbf w \calP_\bbf$ is perfectly placid, i.e., a countable filtered limit of affine pfp schemes whose transition maps are perfections of smooth maps.
\xrema

\theo
\mylabel{convolution.product.witt}
The convolution product $\star$ has the following properties:

\begin{enumerate}
\item
\label{item--convolution.pfp.l-adic}
It is compatible with the convolution product defined in the $\ell$-adic situation \cite[p.~432]{Zhu:Affine} under the $\ell$-adic realization functor $\rho_\ell\co \DM(\calP_\bbf \bsl LG / \calP_\bbf)^\comp \r \D^\bound_{\cstr, \calP_\bbf}(\Fl_\bbf)$ \textup{(}taking values in the bounded derived $\calP_\bbf$-equivariant category of constructible sheaves on $\Fl_\bbf$\textup{)}.

\item 
\label{item--convolution.associative.pfp}
It admits an associativity isomorphism.

\item
\label{item--convolution.DTM}
It preserves the subcategory $\DTM(\calP_\bbf \bsl LG / \calP_\bbf) \subset \DM(\calP_\bbf \bsl LG / \calP_\bbf)$.
For $\calP_\bbf = L^+ G$, it also preserves the abelian subcategory $\MTM(L^+G \bsl LG / L^+G)$.
\end{enumerate}
\xtheo

\pf
Part \refit{convolution.pfp.l-adic} holds since both convolution product functors are defined as the concatenation of the same functors (motivic, resp.~$\ell$-adic) and \myref{DM.perfect.schemes}.\refit{realization}. 
Part \refit{convolution.associative.pfp} results from the base-change isomorphism for proper maps of ind-(pfp schemes) (cf.~\refeq{base.change.pfp}) as in \sref{associative}.

For \refit{convolution.DTM}, we review the key geometric arguments in case $\calP_\bbf$ is the Iwahori group $\calB$, and refer to \sref{convolution.Fl.Tate.2x.Iwahori} for the reduction to this case and for further details concerning the case $\calP_\bbf = \calB$.
Writing $\iota_w : \Fl_w \r \Fl := \Fl_{\bba_0}$ for the inclusion of the stratum, the claim $(\iota_w)_! 1 \star (\iota_{w'})_! 1 \in \DTM(\Fl)$ is first shown in case $w=w'$ is a simple reflection $s$.
In this case the map
\[
\Fl_{\le s} \xtw \Fl_{\le s} \stackrel{(\pr_1, m)} \r \Fl_{\le s} \x \Fl_{\le s} 
\]
is an isomorphism because it is a closed immersion of $2$-dimensional irreducible perfectly proper schemes.
There is a stratification of $\Fl_{\le s} = \Fl_s \sqcup \Fl_e \cong (\A^1)^\perf \sqcup \Spec k$.
Furthermore, the multiplication map $\tilde m : \Fl_{\le s} \xtw \Fl_{\le s} \r \Fl_{\le s}$ is a Whitney--Tate map with respect to the this stratification.
This can be seen as in loc.~cit., except that all the affine lines (intervening as strata) are replaced by their perfections.
The second step uses the isomorphism (see loc.~cit.) $\Fl_w \xtw \Fl_{w'} \r \Fl_{ww'}$ for $w, w' \in W$ such that $l(ww')=l(w)+l(w')$.
This isomorphism is a consequence of \eqref{product_map_eq} and Bruhat--Tits theory, applied to the Witt vector affine flag variety. 
The remaining steps (still for $\calP_\bbf = \calB$) are then formal consequences.

The preservation of $\MTM$ as stated holds because of the corresponding assertion in the $\ell$-adic case \cite[Proposition~2.2]{Zhu:Affine} together with \refit{convolution.pfp.l-adic}.
Here we use \thref{Fl.MTM.adic}.\refit{l-adic.Fl.pfp} to reduce to $\ell$-adic (perverse) sheaves.
\xpf

\section{The motivic Satake equivalence}
We keep the notation from \refsect{Witt.vector.flag.varieties}: $K/\bbQ_p$ is a finite extension with ring of integers $\calO_K$ and finite residue field $k/\bbF_p$.
Let $G$ be a split reductive group over $\calO_K$. 
Fix a maximal split torus contained in a Borel subgroup $T\subset B\subset G$ over $\calO_K$.

Let $\widehat G$ denote the Langlands dual group over $\Q$ formed with respect to the pair $(T,B)$ and a fixed pinning $X\in \on{Lie}(B)$. 
Recall that $(\widehat G,\widehat T, \widehat B, \widehat X)$ is a split reductive $\Q$-group scheme equipped with a pinning whose root datum is dual to the root datum of $(G,T,B)$.
Denote by $\widehat T_\ad$ the image of $\widehat T$ in the adjoint group $\widehat G_\ad$. 
Via $X_*(T)=X^*(\widehat T)$, we can view the half-sum of the $B$-positive roots in $G$ as a cocharacter $\rho\co \GmX \Q \to \widehat T_\ad$.
Then $\GmX \Q$ acts from the right on $\widehat G$ via conjugation as $\on{Ad}_\rho\co (g,\la)\mapsto \rho(\la)^{-1}g\rho(\la)$.
This action preserves the pair $(\widehat T, \widehat B)$. 

\defi 
The {\em extended Langlands dual group} is the $\Q$-group scheme $\widehat G_1:=\widehat G\rtimes^{\on{Ad}_\rho} \GmX \Q$ equipped with the pinning $(\widehat T_1, \widehat B_1, \widehat X_1)$.
\xdefi

Here $\widehat T_1:=\widehat T\rtimes^{\on{Ad}_\rho} \GmX \Q\subset \widehat B\rtimes^{\on{Ad}_\rho}\GmX \Q=:\widehat B_1$ and $\widehat X_1:=(\widehat X, 0)$ is the principal nilpotent element in $\on{Lie}(\widehat B_1)=\on{Lie}(\widehat B)\x \GaX \Q$.
We note that $\widehat G_1$ is a split reductive $\Qq$-group scheme with split maximal torus $\widehat T_1= \widehat T\x \GmX \Q$, since $\on{Ad}_\rho$ acts trivially on $\widehat T$. 
One easily verifies that the map $(g, \la)\mapsto (g\cdot 2\rho(\la),\la^2)$ induces a short exact sequence of $\Qq$-group schemes
\[
1\to \Bmu_2\to \widehat G\x \GmX \Q \to \widehat G_1\to 1,
\]
where $\Bmu_2$ is generated by $(\epsilon, -1)$ with $\epsilon:=(2\rho)(-1)\in \widehat G(\Q)$.
Thus $\widehat G_1$ is the same as the group constructed in \cite[\S2]{FrenkelGross:Rigid} following \cite{Deligne:Letter2007}, cf.~also \cite[Proposition~5.39 ff.]{BuzzardGee:Conjectures} for further examples.
We denote by $d\co \widehat G_1\to \GmX \Q /\Bmu_2\simeq \GmX \Q$ the character where the isomorphism is induced from the square map $\la\mapsto \la^2$.
As an example take $G=\on{PGL}_2$ in which case $\widehat G_1=\on{SL}_2\x^{\Bmu_2}\GmX \Q=\GL_2$, and the character $d$ identifies with the determinant.

We write $\Rep_{\Q}(\widehat G_1)$ for the category of algebraic $\widehat G_{1}$-representations on $\Q$-vector spaces.
The subcategory $\on{Rep}_\Q^{\on{fd}}(\widehat G_1)$ of finite-dimensional representations is a semi-simple abelian category with (absolutely) simple objects
\begin{equation}\label{irred.rep.eq}
V_\mu(n)\defined \Ind_{\widehat{B}_1^\opp}^{\widehat{G}_1}((\mu, n)), \;\;\;\;\mu\in X_*(T)^+, n\in \bbZ, 
\end{equation}
where $\widehat{B}_1^\opp\subset \widehat G_1$ denotes the Borel opposite to $\widehat{B}_1$, and $\mu_n\co \widehat{B}_1^\opp\to \widehat{T}_1\to \GmX \Q$ is the composition of the projection with the character $(\mu,n)\in X_*(T)^+\x \bbZ=X^*(\widehat{T}_1)^+$.
Note that the restriction of $V_\mu(n)$ along $\widehat G\x \GmX \Q \to \widehat G_1$ is the representation $V_\mu\otimes d^{\otimes n}$ where $V_\mu$ is the $\widehat G$-representation of highest weight $\mu$.

Extending scalars from $\Q$ to $\Qq$ we obtain the categories 
\[
\Rep_{\Qq}^{\on{fd}}(\widehat G_1)\subset \Rep_{\Qq}(\widehat G_1)
\]
of algebraic $\widehat G_{1,\Qq}$-representations on $\Qq$-vector spaces which have the same simple objects \eqref{irred.rep.eq}.

The aim of this section is to prove the following theorem due to Zhu \cite{Zhu:Affine} in the context of $\ell$-adic sheaves and in the context of numerical motives \cite[\S2]{Zhu:Geometric}.
It will be stated and proven in a slightly sharper form in \myref{Satake.equiv.thm} and \myref{passing.to.ind.cats} below.

\theo
\mylabel{motivic.satake.witt}
There is an equivalence of $\Qq$-linear abelian tensor categories
\[
\MTM\big(L^+G\bsl LG/L^+G\big)_\Qq \;\overset{\cong}{\lr}\; \Rep_{\Qq}(\widehat{G}_1),
\] 
under which $\IC_\mu(n)\mapsto V_\mu(n)$ for all $\mu\in X_*(T)^+$, $n\in \bbZ$.
For each prime number $\ell\neq p$, this equivalence is compatible under the $\ell$-adic realization with the geometric Satake equivalence constructed in \cite{Zhu:Affine}.
\xtheo

At the left hand side, $\MTM\big(L^+G\bsl LG/L^+G\big)$ denotes the category of stratified mixed Tate motives on the loop group double quotient constructed in \myref{Fl.MTM.adic,DTM.doubly.equivariant}.
The subscript $\Qq$ indicates that we consider this category with $\Qq$-coefficients (as opposed to $\Q$-coefficients like in the preceding sections).
The simple objects $\IC_\mu(n)$, $\mu\in X_*(T)^+$, $n\in \bbZ$ are Tate twists of intersection cohomology motives of the Schubert varieties.
We already know that the convolution product $\star$ for motives on $L^+G\bsl LG/L^+G$ preserves mixed Tate motives (\myref{convolution.product.witt}).
In order to prove the above theorem, we will perform the following steps:
we show in \myref{tannaka.cat} that \begin{equation}\label{neutral.Tannakian.eq}
\big(\MTM\left(L^+G\bsl LG/L^+G\right)^\comp,\star, \omega\big)
\end{equation}
is a neutral Tannakian category. The fiber functor $\omega$ is given by taking global motivic cohomology.
We show that the Tannaka dual group is reductive by showing that the category is semi-simple (\myref{semisimplicity.prop}). This step relies on the Kazhdan--Lusztig parity vanishing and on the semisimplicity of $\MTM(k)$.
In particular, every object in $\MTM\left(L^+G\bsl LG/L^+G\right)$ is a (possibly infinite) direct sum of $\IC_\mu(n)$'s.
We further identify, at least over $\Qq$, the Tannaka dual group $\widetilde G_1 := \Aut^\star(\omega)$ of $\MTM\left(L^+G\bsl LG/L^+G\right)^\comp$ with $\widehat G_1$, by relying on the $\ell$-adic Satake equivalence for the Witt vector affine Grassmannian due to Zhu. 
Once this is done, we obtain \myref{motivic.satake.witt} by passing to ind-completions, cf.~\myref{passing.to.ind.cats}.

\subsection{Semisimplicity} 
\label{sect--semisimplicity}
The Witt vector affine Grassmannian $\Gr_G=(LG/L^+G)^\et$ admits a stratification by left $L^+G$-orbits which yields a presentation as ind-(perfectly proper pfp $k$-schemes) by Witt vector Schubert varieties:
\[
\Gr_G\;=\;\colim_{\mu}\Gr_{G,\leq \mu}.
\]
Here $\mu$ ranges of the partially ordered semigroup $X_*(T)^+$ of $B$-dominant cocharacters.
Recall that $\Gr_{G,\leq \la}\subset \Gr_{G,\leq \mu}$ if and only if $\la \leq \mu$ in the dominance order, i.e., $\mu-\la$ is a sum of positive coroots with non-negative integral coefficients.
The dense open stratum $\Gr_{G, \mu}\subset \Gr_{G,\leq \mu}$ is an irreducible perfectly smooth pfp $k$-scheme of dimension $\lan2\rho,\mu\ran$.
Here $\rho$ denotes the half-sum of the positive roots in $B$ and $\lan\str,\str\ran$ is the natural pairing between the characters and cocharacters.  
In particular, the $L^+G$-orbit stratification satisfies the parity property
\begin{equation}\label{parity.property.eq}
\Gr_{G,\la}\subset \Gr_{G,\leq \mu} \implies \dim( \Gr_{G,\mu})\equiv \dim(\Gr_{G,\la})\mod 2,
\end{equation}
because $\lan2\rho,\mu-\la\ran=2\lan\rho,\mu-\la\ran$ is an even integer.
Recall from \myref{Fl.MTM.adic} the abelian category $\MTM(\Gr_G)$ of mixed stratified Tate motives with respect to the $L^+G$-orbit stratification. 
Its simple objects are the Tate twists of the intersection motives $\IC_\mu(n)$, $\mu\in X_*(T)^+$, $n\in \bbZ$ of the Schubert varieties $\Gr_{G,\leq \mu}$.
By \myref{DTM.doubly.equivariant} the forgetful functor 
\begin{equation}\label{forget.functor.eq}
\MTM(L^+G\bsl LG/L^+G)\;{\lr}\; \MTM(\Gr_G),
\end{equation}
is fully faithful and yields a bijection on simple objects.
The following proposition is analogous to \cite[Proposition~2.1]{MirkovicVilonen:Geometric}. 

\prop
\mylabel{semisimplicity.prop}
The forgetful functor \eqref{forget.functor.eq} is an equivalence of abelian categories.
Both categories are semisimple with simple objects $\IC_\mu(n)$, $\mu\in X_*(T)^+$, $n\in \bbZ$.
\xprop

\pf
The argument follows closely \cite[\S5.1]{RicharzScholbach:Motivic}. We outline the key steps for the reader's convenience. 
If $\MTM(\Gr_G)$ is semisimple, then $\MTM(L^+G\bsl LG/L^+G)$ is semisimple as well by the full faithfulness of \eqref{forget.functor.eq}.
In this case \eqref{forget.functor.eq} is an equivalence because the simple objects agree. 

To prove the semisimplicity of $\MTM(\Gr_G)$ it is enough to show that all extensions between simple objects vanish, i.e., that
\begin{equation}\label{semisimple.eq}
\Ext^1_{\MTM(\Gr_G)}\big(\IC_\la(m), \IC_\mu(n)\big)=0
\end{equation}
for all $\la,\mu\in X_*(T)^+$, $m, n\in \bbZ$.
By the conservativity and t-exactness of the $\ell$-adic realization (\myref{realization.pfp}) and \eqref{parity.property.eq}, the Kazhdan--Lustig parity vanishing holds, i.e., 
\[
{\motH}^i(\iota^*A) = 0,\;\;\;\;\text{whenever \;\; $i\not \equiv 0\mod 2$},
\]
where $A\in \MTM(\Gr_G)$ and $\iota\co X\subset \Gr_G$ is a finite union of Schubert varieties and where ${\motH}^i$ denotes the truncation with respect to the motivic $t$-structure on $X$.
This formally implies 
\[
\Ext^1_{\MTM(\Gr_G)}\big(\IC_\la(m), \IC_\mu(n)\big) = \left\{
\begin{array}{ll}
\Ext_{\MTM(k)}(1_k(m), 1_k(n)) & \text{ if $\la = \mu$} \\
0 & \, \textrm{else} \\
\end{array}
\right. 
\]
where we have used that $\MTM(\Gr_\mu)=\MTM(k)$, cf.~the argument in \sref{Hom.IC}.
Finally, the higher algebraic $K$-groups of the finite field $k$ are torsion (i.e., their rationalizations vanish), so that $\MTM(k)$ is semisimple.
This shows \eqref{semisimple.eq}.
\xpf

\subsection{The tensor structure} 
\label{sect--tensor}

By \myref{convolution.product.witt} the category $\MTM(L^+G\bsl LG/ L^+G)$ is stable under the convolution product $\star$ compatible with the $\ell$-adic realization, i.e., for each pair of such motives $A, B$, there is a functorial isomorphism of $\ell$-adic perverse sheaves
\begin{equation}\label{ell.adic.convolution}
\rho_\ell(A\star B) \cong \rho_\ell(A)\star_\ell \rho_\ell(B),
\end{equation}
where 
$$\rho_\ell\co \MTM(L^+ G \bsl LG / L^+ G)^\comp \r \Perv_{L^+ G}(\Gr_G, \Ql)\eqlabel{rho.ell.mal.wieder}$$
denotes the $\ell$-adic realization taking values in the category of $L^+G$-equivariant perverse $\ell$-adic sheaves on $\Gr_G$ (cf.~\myref{convolution.product.witt}), and where $\star_\ell$ denotes the convolution of $\ell$-adic complexes. 

\prop 
\mylabel{constraints.prop}
The convolution product $\star$ can be upgraded to a unique symmetric monoidal structure on $\MTM(L^+ G \bsl LG / L^+ G)$ such that the $\ell$-adic realization functor $$\rho_\ell\co \MTM(L^+ G \bsl LG / L^+ G)^\comp \r \Perv_{L^+ G}(\Gr_G, \Ql)$$ is a symmetric monoidal functor, where the target category carries the symmetric monoidal structure established in the geometric Satake equivalence \cite[Proposition~2.21]{Zhu:Affine}.
\xprop
\pf
There exist functorial commutativity and associativity constraints 
\[
c_{A,B}\co A\star B\cong B\star A, \;\;\;\text{and}\;\;\; a_{A,B,C}\co (A\star B)\star C\cong A\star(B\star C),
\]
which are uniquely determined by the following two properties:
\begin{enumerate}
\item[i\textup{)}] The isomorphisms are colimit-preserving in each argument.
\item[ii\textup{)}] The constraints map under the $\ell$-adic realization \refeq{rho.ell.mal.wieder} to the constraints used in geometric Satake as in \cite[Proposition~2.21]{Zhu:Affine}.
\end{enumerate}

The uniqueness follows from the faithfulness of the $\ell$-adic realization (\myref{realization.pfp}).
The associativity constraint is constructed in \myref{convolution.product.witt}.
The construction of the (correct) commutativity constraint is the most subtle part and follows \cite[\S2.4.3]{Zhu:Affine}.
For details the reader is referred to \sref{constraints}.
Here we only sketch the main ideas. 
The anti-involution $\theta\co G\to G$ is defined by $g\mapsto (g^*)^{-1}=(g^{-1})^*$, where $(\str)^*$ denotes the Cartan involution.
By functoriality, we obtain an anti-involution on $LG$ preserving $L^+G$, and thus an equivalence of prestacks, still denoted
\[
\theta\co L^+G\bsl LG/L^+G \stackrel {\simeq}{\longrightarrow} L^+G\bsl LG/L^+G.
\]
Then for $A,B\in \MTM(L^+G\bsl LG/L^+G)$ we have a functorial isomorphism $\theta^!(A\star B)\cong (\theta^!B)\star (\theta^!A)$.
Also we have a (carefully chosen) equivalence $\theta^!\simeq \id$ of endofunctors on $\MTM(L^+G\bsl LG/L^+G)$.
The point being that $\theta$ preserves Schubert varieties and hence also intersection motives. 
However, to obtain the correct equivalence $\theta^!\simeq \id$ compatible with the $\ell$-adic realization one needs to introduce a carefully chosen sign, cf.~the aforementioned references.  
Combining the above equivalences leads to the commutativity constraint 
$c'_{A, B}\co A\star B\cong B\star A$.
We finally define
\[
c_{\IC_\mu(n)\star\IC_\la(m)}:=(-1)^{\lan2\rho, \mu +\la\ran}c'_{\IC_\mu(n)\star\IC_\la(m)},
\] 
and extend it linearly to all objects in $\MTM(L^+G\bsl LG/L^+G)$ using the semisimplicity, cf.~\myref{semisimplicity.prop}. 
\xpf

\subsection{The Tannakian structure} 
\label{sect--tannakian}
Next we equip the $\bbQ$-linear, abelian, symmetric monoidal, semisimple category $(\MTM(L^+G\bsl LG/L^+G), \star)$ with a tensor functor valued in vector spaces $(\Vect_\bbQ,\otimes)$ leading to the aforementioned Tannakian structure.
Following \cite[\S5.3]{RicharzScholbach:Motivic} we define the {\it fibre functor} 
as the composition
$$
\omega\co \MTM(L^+G\bsl LG/L^+G) \cong  \MTM (\Gr_G) \stackrel{\epsilon_!} \lr \DTM(k) \stackrel{\gr^\cl} \lr \MTM(k) \stackrel{\gr^\W} \lr \MTM(k)^{\wt = 0} \cong \Vect_\Q
$$
of the forgetful functor (which is an equivalence by \myref{semisimplicity.prop}; note that mixed Tate motives are with respect to the stratification by $L^+ G$-orbits), the pushforward along the structural map $\epsilon\co \Gr_G \r \Spec(k)$
(which preserves Tate motives by \inter{Lemma~\iref{lemm--Tate.proper.descent}}, using that the stratification of $\Gr_G$ by $L^+G$-orbits is perfectly cellular), followed by the grading functors for the classical (which agrees in this case with the motivic) $t$-structure $\gr^\cl$ and the weight structure $\gr^\W$, and finally the equivalence \cite{Levine:Tate} of pure Tate motives of weight 0 with the category of $\Q$-vector spaces.

\lemm
\mylabel{monoidal.lemm}
The fibre functor $\omega \colon \MTM(L^+G\bsl LG/L^+G)\to \Vect_\bbQ$ has a natural monoidal structure with respect to the convolution product on the source and the ordinary tensor product on the target which is compatible with the $\ell$-adic realization. 
\xlemm
\pf
At least on the level of homotopy categories the composition of functors
\[
\DM(L^+G\bsl LG/L^+G)\lr \DM(\Gr_G)\stackrel {\epsilon_!}\lr \DM(k)
\]
has a natural monoidal structure with respect to convolution on the source and the tensor product on the target, cf.~\sref{pushforward.symmetric.monoidal} whose arguments translate literally. 
Next the conservativity of the $\ell$-adic realization implies the Kazdhan--Lusztig parity vanishing ${^\cl\H^i}(\epsilon_!\IC_\mu(n))=0$ whenever $i\not \equiv \lan2\rho,\mu\ran \mod 2$.
Hence, the composition 
\[
\MTM(L^+G\bsl LG/L^+G) \cong  \MTM (\Gr_G) \stackrel{\epsilon_!} \lr \DTM(k) \stackrel{\gr^\cl} \lr \MTM(k)
\]
admits a natural monoidal structure as well. 
The vanishing is needed in order to obtain a monoidal structure as opposed to a $\bbZ/2$-graded monoidal structure coming from a sign in the formation of tensor products of complexes. 
Using the natural monoidal structure on the weight graduation functor $\gr^W$ we obtain a monoidal structure on $\omega$.
The compatibility with the $\ell$-adic realization follows from \myref{DM.perfect.schemes}.\refit{realization} and the fact that the $\ell$-adic fiber functor is defined by the analogous functors.
\xpf

\theo
\mylabel{tannaka.cat}
The category $\MTM(L^+G\bsl LG/L^+G)^\comp$ of compact mixed Tate motives on the double quotient, endowed with the convolution product, the constraints from \myref{constraints.prop} and $\omega$ as fibre functor
is a neutral Tannakian category over $\bbQ$ \textup{(}\cite[Chapter~II, Definition~2.19]{DeligneMilneOgusShih:Hodge}\textup{)}.
\xtheo
\pf
Note that the restriction of $\omega$ to compact objects takes values in finite-dimensional vector spaces.
As in \sref{Tannaka_Cat} we check the conditions in \cite[Chapter~II, Proposition~1.20]{DeligneMilneOgusShih:Hodge}:
\begin{enumerate}
\item {\it The functor $\omega$ has the structure of a tensor functor.} This is \myref{monoidal.lemm}.

\item {\it The functor $\omega$ is $\bbQ$-linear, exact and faithful.} The functor $\omega$ is clearly $\Q$-linear and additive. 
Hence, it is exact because $\MTM(L^+G\bsl LG/L^+G)^\comp$ is semi-simple. 
The faithfulness follows from the compatibility with the $\ell$-adic realization, the faithfulness of the latter, and the faithfulness of the $\ell$-adic fiber functor $\omega_\ell$ \cite[Corollary~2.10]{Zhu:Affine}.

\item {\it The constraints constructed in \myref{constraints.prop} give the usual constraints in $\Vect_\Q$ after applying $\omega$.} This is immediate from \myref{monoidal.lemm} and the $\ell$-adic case, cf.~\cite[Proposition~2.21]{Zhu:Affine}.

\item {\it Neutral object.} This is the skyscraper $\IC_0$ supported at the base point.

\item {\it Any $M\in \MTM(L^+G\bsl LG/L^+G)^\comp$ with $\dim_\bbQ \omega(M)=1$ admits a dual object $M^{-1}$ such that $M \star M^{-1} = \IC_0$.} 
Any such object is necessarily of the form $M=\IC_\mu$ with $0$-dimensional support, in which case $M^{-1}=\IC_{-\mu}$.

\end{enumerate}
\xpf

\subsection{The Tannakian group}

\label{sect--dual.group}

\theo
\mylabel{Satake.equiv.thm}
There is an equivalence of Tannakian categories
\[
\left (\MTM(L^+G\bsl LG/L^+G)_\Qq^\comp,\;\star,\;\omega\right ) \;\simeq\; \left(\Rep_{\Qq}^{\on{fd}}(\widehat{G}_1),\;\t, \;v\right),
\]
where the subscript $\Qq$ indicates the category with the same objects, but $\Hom$-spaces are tensored with $\Qq$.
At the right, $v\co \Rep_{\Qq}^{\on{fd}}(\widehat{G}_1)\to \on{Vec}_\Qq$ denotes the forgetful functor.
The intersection motives $\IC_\mu(n)$ correspond to the simple $\widehat{G}_1$-representations $V_\mu(n)$ for $(\mu,n)\in X_*(T)^+\x \bbZ=X^*(\widehat{T}_1)^+$.
\xtheo
\pf 
The argument is very similar to the proof of \sref{Satake}; it relies on the $\ell$-adic geometric Satake equivalence proven in \cite{Zhu:Affine}.
Let us recall the key points. 
Denote by $\widetilde G_1:=\Aut^\star(\omega)$ the affine $\Qq$-group scheme of tensor automorphisms provided by \myref{tannaka.cat}. 
We need to prove an isomorphism $\widetilde G_1\simeq \widehat G_1$ of $\Qq$-group schemes. 

For a prime number $\ell\nmid p$, the $\ell$-adic realization induces an equivalence of Tannakian categories
\begin{equation}
\label{equiv.realization.eq}
\big(\MTM(L^+G\bsl LG/L^+G)_\Qq^\comp, \star, \omega\big)\otimes_\Qq \Qq_\ell\; \stackrel\cong\lr \;\big(\on{Sat}_{G,\ell}, \star_\ell, \omega_\ell\big),
\end{equation}
where $\on{Sat}_{G,\ell}$ denotes the full semi-simple subcategory of $L^+G$-equivariant $\Qq_\ell$-adic perverse sheaves sheaves on $\Gr_G$ generated by the Tate twisted intersection complexes $\IC_{\mu,\ell}(n)$, $\mu\in X_*(T)^+$, $n\in \bbZ$. 
The affine $\Qq$-group scheme of tensor automorphisms $\Aut^\star(\omega_\ell)$ is isomorphic to $\widehat G_1\otimes \Qq_\ell$ by the $\ell$-adic geometric Satake equivalence as in \cite{Zhu:Affine}. 
This reference is for the ordinary affine Grassmannian, but the arguments carry over to the case of Witt vectors, see also \cite[\S2]{Zhu:Geometric}. 
Hence, \eqref{equiv.realization.eq} implies that there is an isomorphism of $\Qq_\ell$-group schemes
\[
\widetilde G_1\otimes \Qq_\ell \;\simeq \; \widehat G_1\otimes \Qq_\ell
\]
for all prime numbers $\ell\nmid p$. 
In particular, $\widetilde G_1$ is a reductive $\Qq$-group scheme by fpqc descent along the extension $\Qq_\ell/\Qq$ for one fixed prime number $\ell\nmid p$.
We conclude $\widetilde G_1\simeq \widehat G_1$ by the isomorphism theorem for split reductive groups.
\xpf

\rema
\mylabel{classical.satake.rema}
The relation with the classical Satake isomorphism \cite{Gross:Satake} is as follows, cf.~also \cite[\S6.4]{RicharzScholbach:Motivic} for a detailed discussion in the (completely analogous) equal characteristic case.
Consider the spherical Hecke ring
\[
\calH_{G}\defined \calC^c\big(G(K)/\!/G(\calO_K); \bbZ\big)
\]
of functions $G(\calO_K)\bsl G(K)/G(\calO_K)\to \bbZ$ supported on finitely many double cosets with ring structure given by the convolution product. 
Taking the trace of geometric Frobenius on motives as in \cite{Cisinski:SurveyCoho} induces a surjective ring homomorphism from the Grothendieck ring
\[
K_0 \left (\MTM(L^+G\bsl LG/L^+G)^\comp \right )\;\to\; \calH_G\otimes_\bbZ \bbZ[q^{-1}], \; M\mapsto f_M,
\]
with kernel generated by the class $[\IC_0(-1)]-q[\IC_0]$ (the trace of geometric Frobenius on $\Qq(-1)$ is given by multiplication with $q$).
Hence, under \myref{Satake.equiv.thm} the map $[\IC_\mu(n)]\mapsto [V_\mu(n)]$ induces an isomorphism of rings
\[
\calH_{G}\otimes \bbZ[q^{-1}] \stackrel \cong \lr R(\widehat G_1)/([d^{-1}]-q),
\]
where $R(\widehat G_1)=K_0\Rep_{\Qq}^{\on{fd}}(\widehat{G}_1)$. 
This is immediate from the preceding discussion using $[V_0(-1)]=[d^{-1}]$ in $R(\widehat G_1)$. 
Now choosing a square root $q^{1/2}\in \Qq$ the composition gives the classical Satake isomorphism as in \cite[Proposition~3.6, (3.12)]{Gross:Satake},
\[
\calH_{G}\otimes \bbZ[q^{\pm1/2}] \cong \big(R(\widehat G_1)/([d^{-1}]-q)\big)[q^{\pm 1/2}] \cong R(\widehat G)\otimes \bbZ[q^{\pm 1/2}],
\]
where the last isomorphism is induced from $[V_\mu(0)]\mapsto q^{\lan\rho,\mu\ran}[V_\mu]$.
\xrema

\rema
\mylabel{passing.to.ind.cats}
\myref{Satake.equiv.thm} implies \myref{motivic.satake.witt} by passing to ind-completions.
Indeed, the category of algebraic representations $\Rep_\Qq(\widehat G_1)$ is compactly generated by its subcategory $\Rep_\Qq^{\on{fd}}(\widehat G_1)$ of finite-dimensional representations.
Similarly, $\MTM(L^+G\bsl LG/L^+G)$ is compactly generated by $\MTM(L^+G\bsl LG/L^+G)^\comp$.
\xrema

\rema
\mylabel{extension.witt.satake}
We expect that the analogue of \myref{Satake.equiv.thm} for $\Q$-coefficients (as opposed to $\Qq$-coefficients) holds as well. 
For this one needs to show that the Tannakian group $\widetilde G_1=\Aut^\star(\omega)$ is $\Q$-split.
It will then be automatically isomorphic to $\widehat G_1$ by the isomorphism theorem. 
In order to obtain a more canonical identification $\widetilde G_1\simeq \widehat G_1$ it would be interesting to construct constant term functors in this setting. 

Also, we expect that \myref{Satake.equiv.thm} admits, similarly to \cite{Zhu:Affine}, a generalization from split reductive groups to unramified reductive groups.  
For this a key input is the existence of motivic t-structures on categories of stratified Artin--Tate motives (as opposed to Tate motives). 
Here we expect the methods of \inter{\S 3} to carry over. 
\xrema

\bibliographystyle{alphaurl}
\bibliography{bib}

\end{document}